\documentclass[a4paper,11pt]{amsart}
\usepackage{mathrsfs}
\usepackage{appendix}  
\usepackage{amssymb}  
\usepackage{fancyhdr}      
\usepackage{charter} 
\usepackage{typearea}

\usepackage[a4paper,top=2cm,bottom=2.5cm,left=2cm,right=2cm]{geometry}

\usepackage{amsmath,amsfonts,amscd,amsthm}
\usepackage{color}                    
\usepackage[colorlinks,linkcolor=blue,anchorcolor=blue,citecolor=blue]{hyperref} %for show the ref

\newcommand{\dbar}{\ensuremath{\overline\partial}}
\newcommand{\dbarstar}{\ensuremath{\overline\partial^*}}
\newcommand{\C}{\ensuremath{\mathbb{C}}}

\newcommand{\R}{\ensuremath{\mathbb{R}}}

\newcommand{\ov}{\overline}
 
\newcommand{\ddbar}{\overline\partial}
\newcommand{\til}[1]{\widetilde{#1}}

\newcommand{\norm}[1]{\left\Vert#1\right\Vert}
\newcommand{\abs}[1]{\left\vert#1\right\vert}

\newcommand{\set}[1]{\left\{#1\right\}}
\newcommand{\To}{\rightarrow}

\newcommand{\cali}[1]{\mathscr{#1}}

\newcommand{\be}{\begin{eqnarray}}
\newcommand{\ee}{\end{eqnarray}}

\newcommand{\comment}[1]{}

\DeclareFontFamily{U}{mathx}{\hyphenchar\font45}
\DeclareFontShape{U}{mathx}{m}{n}{
	<5> <6> <7> <8> <9> <10>
	<10.95> <12> <14.4> <17.28> <20.74> <24.88>
	mathx10
}{}
\DeclareSymbolFont{mathx}{U}{mathx}{m}{n}
\DeclareFontSubstitution{U}{mathx}{m}{n}
\DeclareMathAccent{\widecheck}{0}{mathx}{"71}
\DeclareMathAccent{\wideparen}{0}{mathx}{"75}

\def\omz{\Omega}

\makeatletter
\newcommand{\sumprime}{\if@display\sideset{}{'}\sum%
	\else\sum'\fi}
\makeatother

\newtheorem{thm}{Theorem}[section]
\newtheorem{prop}[thm]{Proposition} 
\newtheorem{lem}[thm]{Lemma}
\newtheorem{cor}[thm]{Corollary}

\theoremstyle{definition}

\theoremstyle{remark}
\newtheorem{rem}[thm]{\bf Remark}
\numberwithin{equation}{section}

\providecommand\ufootnote[1]{{\let\thefootnote\relax\footnote[0]{#1}}}

\newcommand{\fc}{\mathcal F}

\newcommand{\bs}{\mathbb S}
\newcommand{\bt}{\mathbb T}

\newcommand{\cO}{\mathscr{O}}

\newcommand{\N}{\mathbb{N}}

\newcommand{\ol}{\overline}

\newcommand{\pa}{\partial}

\newcommand{\wt}{\widetilde}

 \DeclareMathOperator{\Dom}{Dom}
  
\DeclareMathOperator{\Tr}{Tr} 
\DeclareMathOperator{\vol}{vol}

\begin{document} 
\title{
Heat kernel asymptotics for Kodaira Laplacians of high power of line bundle over complex manifolds
   }
	\author{Huan Wang} 
	\address{Huan Wang, Institute of Mathematics, Henan Academy of Sciences,
		Zhengzhou 450046, Henan, China}   
	\email{huanwang@hnas.ac.cn, huanwang2016@hotmail.com} \thanks{The first author is supported by the High-level Talent Research Start-up Project Funding of the Henan Academy of Sciences (Project No. 241819110) and the Guest Program of MPIM Bonn. The second author is partially supported by the Austrian Science Fund (FWF) grant No. 10.55776/ESP367.}     
	\author{Weixia Zhu}   
	\address{Weixia Zhu, Faculty of Mathematics, University of Vienna, 
		Vienna 1090, Austria}   
	\email{weixia.zhu@univie.ac.at, zhuvixia@gmail.com}
 
 \date{2024.Nov.17}
 
\begin{abstract}   
This paper presents a simple method to prove the heat kernel asymptotics for the Kodaira Laplacian with respect to the high power of a holomorphic Hermitian line bundle $(L,h^L)$ over a possibly non-compact Hermitian manifold $(M,\omega)$. As a consequence, we give a direct proof of the holomorphic Morse inequalities on covering manifolds. Furthermore, we generalize it to the vector bundle via the $L^2$ Le Potier isomorphism and provide an algebraic version of the holomorphic Morse inequalities. The approach used in this work employs a scaling technique and is applicable to $M$ regardless of its compactness. 
 
		
\smallskip
		
\end{abstract} 

\maketitle  	 
\tableofcontents     
	 
\section{Introduction}
     
The asymptotic expansion of the heat kernel provides a unified framework for various topics, including the Atiyah–Singer index theorem, holomorphic Morse inequalities, Bergman kernel asymptotics, and analytic torsion \cite{Ma20}. Let $(M, \omega)$ be a complex manifold of dimension $n$ with Hermitian metric $\omega$. Consider a holomorphic Hermitian line bundle $(L, h^L)$ on $M$, and let $\Box^q_{k}$ denote the Gaffney extension of the Kodaira Laplacian acting on $(0,q)$-forms with values in the $k$-th tensor power of $L$, $0 \leq q \leq n$. For $t>0$, let $e^{-\frac{t}{k} \Box^q_{k}}(z,z')$ be the heat kernel associated with the rescaled Kodaira Laplacian $\frac{1}{k}\Box^q_{k}$ on $M$. The asymptotic behavior of this heat kernel as $k \to +\infty$ has been demonstrated by Bismut through probability theory, by Demailly using Mehler’s formula and asymptotic estimates, and by Ma-Marinescu via Bismut-Lebeau’s analytic localization techniques. %These results hold significant importance in various areas, including holomorphic Morse inequalities, asymptotic expansions of the Bergman kernel and the asymptotic of the analytic torsion.

The primary purpose of this paper  
 is to provide an alternative proof of this classical result \cite[Theorem 1.5]{B87} \cite[Theorem 1.6.1]{MM07}, via the scaling technique motivated by the establishment of heat kernel asymptotics for the Kohn Laplacian on CR manifolds \cite{HZ23} (See also \cite{B04,HM12}). This technique can also be applied to compact manifolds with boundary \cite{HMZ24}. Let $\nabla^L$ denote the Chern connection of {the line bundle} $(L, h^L)$ and let $R^L=(\nabla^L)^2$ be its curvature. Moreover, let $\dot{R^L} \in \text{End}(T^{1,0}M)$ and $\varTheta$ be given by
\[
\langle \dot{R^L} u,v \rangle_\omega = R^L(u, \overline{v}),\quad \varTheta = -\sum_{i,j=1}^n R^L(Z_i, \overline{Z}_j) \overline{w}^i \wedge (\overline{w}^j \wedge)^*,
\]
where $\{Z_j\}_{j=1}^n$ is a local orthonormal frame of $T^{(1,0)}M$, $\{w^j\}_{j=1}^n$ is the dual frame, and $(\overline{w}^j \wedge)^* = i_{\overline{Z}_j}$ is the contraction.
 Our main theorem can be formulated as follows:

\begin{thm}\label{heatkernelthm}
Let $(M,\omega)$ be a (not necessarily compact) Hermitian manifold of dimension $n$. Let $(L,h^L)$ be a holomorphic Hermitian line bundle on $M$. 
Let $I\subset\mathbb \R_+=(0,\infty)$  be a compact interval and let $K\subset M$ be a compact set. Then, there is a constant $C>0$ independent of $k$ such that, for all $z\in K$, $t\in I$, $0\leq q\leq n$, we have
\begin{equation}\label{eq_ub}
\abs{k^{-n}e^{-\frac{t}{k} \Box^q_{k}}(z,z)}_{\mathscr L(T^{*0,q}_zM,T^{*0,q}_zM)}\leq C,
\end{equation} 
{where $|\cdot|_{\mathscr L}$ is defined in \eqref{eq_L}.}
Moreover, at each point $z\in M$, we have 
\begin{equation}\label{eq_hka}
\lim_{k\to+\infty}k^{-n}e^{-\frac{t}{k}\Box^q_{k}}(z,z)=\frac{\det(\dot{R^L}/2\pi)\exp(t\varTheta)}{\det(1-\exp(-t\dot{R^L}))}(z).
\end{equation}    		
\end{thm}
 When $M$ is compact, Theorem \ref{heatkernelthm} reduces to the classical result \cite[Theorem 1.5]{B87} \cite[Theorem 1.6.1]{MM07}.
When $M$ is a relatively compact domain, Theorem \ref{heatkernelthm} corresponds to \cite[(3.2.32)]{MM07} and plays a crucial role in proving the abstract Morse inequalities for the $L^2$-cohomology.
Our proof involves three ingredients: the heat kernel in model case \eqref{e-gue210523yydII}, the localization process \eqref{eq_scaling} within coordinate neighborhoods \eqref{local_co}, and elliptic estimates for the scaled Kodaira Laplacian \eqref{boxrhok}. Ma-Marinescu's localization technique \cite[\S 1.6]{MM07} relies on establishing a normal coordinate chart through geodesics and parallel transport. Berndttson-Berman provide an relatively elementary localization method for Bergman kernels \cite{B04}. In our proof, we use for the first time the scaling techniques on the asymptotics of heat kernel over complex manifolds.

\begin{rem}
In Theorem \ref{heatkernelthm}, if an eigenvalue of $\dot{R^L}$ at $p\in M$ is zero, then its contribution to the term $\det(\dot{R^L}/2\pi)/\det(1 -\exp(-t \dot{R^L}))$ is $1/(2\pi t)$. Specifically, if all eigenvalues vanish at $p$, we obtain
\[
\lim_{k\to+\infty}k^{-n}e^{-\frac{t}{k}\Box^q_{k}}(p,p) = 1/(2\pi t)^n.
\]
A key advantage of our proof is that it not only recovers the above case but also provides more information for the case when $R=0$ at $p$ (see Sec.\,\ref{sec_rl0} for details). If the vanishing order $\sigma$ of the local weight is known, by adjusting the scaling map, we can obtain finer estimate for the coefficient of $k^m$ in the asymptotic expansion for $m < n$. More precisely, if near $p = 0$, the local weight $\phi(z)$ of $h^L$ can be expressed as $\phi(z) = \sum_{j=1}^n \lambda_j |z_j|^\sigma + O(|z|^{\sigma+1})$, where $\sigma \geq 3$. From our approach, by constructing a new scaling map \eqref{aaaa3}, we derive the following result:
\[
\lim_{k \to \infty} \frac{1}{k^{2n/\sigma}} e^{-\frac{t}{k^{2/\sigma}} \Box_k^q}(p, p) = e^{-t \Box_{\phi_0}^q}(0,0),
\]
where $e^{-t \Box_{\phi_0}^q}(z,z')$ is the heat kernel of the weighted Laplacian $\Box_{\phi_0}^q$ in $\mathbb{C}^n$, with $\phi_0 =\sum_{j=1}^n\lambda_j|z_j|^\sigma$. It is a challenge problem to calculate the heat kernel for $\Box_{\phi_0}^q$. However, when the time parameter $t\to +\infty$, its heat kernel converges to be the Bergman kernel $B_{\phi_0}^{q}$, which is computable in certain cases (such as \cite{MS23}). Our method provides a simple yet effective way to study the heat kernel, even at degenerate points.

\end{rem}

In addition to providing a new proof of the heat kernel asymptotics, Theorem \ref{heatkernelthm} has many important applications on both compact and non-compact manifolds. For instance, it can be used to derive Demailly's holomorphic Morse inequalities on compact complex manifolds \cite{Dem:85}, and, in fact, it provides a broader version of these inequalities. By applying the theorem to a fundamental domain of a covering manifold and using the Lebesgue dominated convergence theorem, we obtain immediately the holomorphic Morse inequalities on covering manifolds of Chiose-Marinescu-Todor (see Theorem \ref{thm_cplxcover}). Additionally, by establishing the $L^2$ version of the Le Potier isomorphism (see Proposition \ref{eq_l2lep}), we derive a new holomorphic Morse inequality in the following:

\begin{thm}\label{cor_sym} 
Let $(\wt M,\widetilde \omega)$ be a Hermitian manifold of dimension $n$. Suppose that $\Gamma$ is a discrete group acting  holomorphically, freely, and properly on $\widetilde M$, such that $\widetilde{\omega}$ is a $\Gamma$-equivariant Hermitian
metric and the quotient $M=\widetilde M/\Gamma$ is compact. Let $(\widetilde E,h^{\widetilde E})$ be a $\Gamma$-equivariant
holomorphic Hermitian vector bundle of rank $r$ on $\widetilde M$. Then, as $k\rightarrow \infty$,  the von Neumann dimension of the space of $L^2$ holomorphic sections satisfies 
\begin{equation*}%\label{eq_symhol}  
\dim_{\Gamma} H_{(2)}^{0} (\til M, S^k(\til E))\geq  \frac{k^{(n+r-1)}}{(n+r-1)!}\int_{P(E^*)(\leq 1)} c_1(\cO_{P(E^*)}(1),h^{\cO_{P(E^*)}(1)})^{n+r-1}+o(k^{n+r-1}). 
\end{equation*}  
\end{thm}
Also, from Theorem \ref{thm_cplxcover}, we establish a equality between the volume of a $\Gamma$-nef line bundle $\wt L$ and the volume of $L=\wt L\setminus \Gamma$ on K\"{a}hler coverings: 
\begin{thm}\label{thm_cover}
If $\til M$ is a K\"{a}hler covering manifold of dimension $n$ and $\til L$ is $\Gamma$-nef, then 
 \begin{equation*}
 \vol_{\Gamma,(2)}(\til L)=\vol(L)=\int_M c_1(L)^n.
 \end{equation*}	
 Moreover, if $\til L_1$ and $\til L_2$ are $\Gamma$-nef, then $\vol_{\Gamma,(2)}(\til L_1\otimes \til L_2)\geq \vol_{\Gamma,(2)}(\til L_1)+\vol_{\Gamma,(2)}(\til L_2)$.
\end{thm} 
When $\Gamma$ is trivial (i.e., $\widetilde M=M$), Theorem \ref{cor_sym} reduces to Demailly \cite[7]{Dem:91}, and Theorem \ref{thm_cover} reduces to \cite[Corollary 2.3.38]{MM07}. When $r=1$ (i.e., $E$ is a line bundle), Theorem \ref{cor_sym} was estsblished in \cite[Corollary 3.6.8]{MM07}.

Our paper is structured as follows: In Sec.\,\ref{sec:prelim}, we revisit the necessary definitions and establish basic settings. In Sec.\,\ref{estimate}, we provide a proof for Theorem \ref{heatkernelthm}. To achieve this, we introduce a scaled Laplacian denoted as $\Box^q_{(k)}$ and its corresponding distribution kernel $A_{(k)}(t,z,z')$ in Sec.\,\ref{scas}. We proceed to demonstrate that $A_{(k)}(t,z,z')$ is uniformly bounded in $k$ within the $\mathscr{C}^\infty$ topology, and it uniformly converges to the heat kernel of certain deformed Laplacian in $\mathbb{C}^n$. In Sec.\,\ref{sec-morse}, we prove Theorem \ref{thm_cplxcover} and its CR counterpart, Theorem \ref{cor_sym}, Theorem \ref{thm_cover}, and algebraic Morse inequalities on K\"{a}hler coverings, see Theorm \ref{thm_amikc}.
    
\section{Basic setup}\label{sec:prelim} 
	     
Let $M$ be a complex manifold of dimension $n$, equipped with a Hermitian form $\omega$. The corresponding volume form is given by $dv_M = \omega^n/n!$. Denoted by $\langle\,\cdot\, , \,\cdot\,\rangle = \langle\,\cdot \, ,\,\cdot\,\rangle_\omega$ the induced Hermitian metric on the holomorphic tangent bundle $T^{1,0}M$. This metric $\langle\,\cdot\, , \,\cdot\,\rangle$ is also used to define the induced Hermitian metric on the space $\Lambda^{p,q}T^*M$ of $(p,q)$-forms on $M$. Let $\mathscr L(V, V)$ be the space of linear endomorphism on a vector space $V$.
For $f(z)\in \mathscr L(T^{*0,q}_zM,T^{*0,q}_zM)$, set
\begin{equation}\label{eq_L}
\abs{f(z)}_{\mathscr L(T^{*0,q}_zM,T^{*0,q}_zM)}:=\sumprime_{J,K}\abs{\langle\,f(z)\ol 
 w^J(z) ,  \ol w^K(z)\,\rangle},
\end{equation} 
and the trace as $\Tr_q f(z):=\sumprime_{J}\abs{\langle\,f(z)\ol w^J(z) ,  \ol w^J(z)\,\rangle}$,
where $\set{\ol w^J(z)}_{|J|=q}$ is an orthonormal basis of  $T^{*0,q}_zM$ with respect to $\langle\,\cdot\, , \,\cdot\,\rangle$. 
	 
Let $D\subset M$ be an open subset. We define $\omz^{0,q}(D)$ and $\omz_c^{0,q}(D)$, $0\leq q\leq n$, as the space of smooth sections of $T^{*0,q}M$ over $D$ and the subspace of $\omz^{0,q}(D)$ whose elements have compact support in $D$, respectively.  The induced $L^2$-inner product is denoted by $(\,\cdot\, ,\,\cdot\,)=\int_D \langle\,\cdot\, , \,\cdot\,\rangle dv_M$. Define $L^2_{(0,q)}(M)$ as the completion of $\omz_c^{0,q}(M)$ with respect to $(\,\cdot\, ,\,\cdot\,)$ and $\|\cdot\|$ as the $L^2$-norm on $L^2_{(0,q)}(M)$.

Consider a holomorphic Hermitian line bundle $(L,h^L)$ over $M$, where $h^L$ denotes a Hermitian metric of $L$. Let $s$ be a local holomorphic frame of $L$ on $D$ such that $|s(z)|^2_{h^L}=e^{-\phi(z)}$ for $z\in D$, where $\phi\in\mathscr C^\infty(D,\mathbb R)$. For $k\in\N\setminus \{0\}$, let $L^k=L^{\otimes k}$ be the $k$-th tensor power of $L$ over M. Denoted by $h^{L^k}$ the induced Hermitian metric on $L^k$ and $s^k$ the induced local holomorphic frame of $L^k$ on $D$. We define $\omz^{0,q}(D,L^k)$ as the space of smooth sections of $T^{*0,q}M\otimes L^k$ over $D$, and $\omz_c^{0,q}(D,L^k)$ as the subspace of $\omz^{0,q}(D,L^k)$ whose elements have compact support in $D$. For $u\in\omz^{0,q}(D,L^k)$, we write $|u|:=|u(z)|^2_{h^{L^k}}$ for simplicity. Let $$\dbar_{k} \colon \omz^{0,q}(M,L^k)\to \omz^{0,q+1}(M,L^k)$$
be the Cauchy-Riemann operator acting on $(0,q)$-forms with values in $L^k$ such that, locally,
$$
\dbar_{k}(s^ku)=s^k\dbar u
$$ 
on $D$ for $s^ku:=s^k\otimes u\in\omz^{0,q}(D,L^k)$ with $u\in \omz^{0,q}(D)$. Here  $\dbar$ is the Cauchy-Riemann operator acts on $\omz^{0,q}(M)$.

We define the $L^2$-inner product on $\omz^{0,q}(M,L^k)$ as $(\,\cdot\, ,\,\cdot\,)_{{k}}$. Let $L^2_{(0,q)}(M,L^k)$ denote the completion of $\omz_c^{0,q}(M,L^k)$ with respect to $(\,\cdot\, ,\,\cdot\,)_{{k}}$ and $\|\cdot\|_{{k}}$ the $L^2$-norms on  $L^2_{(0,q)}(M,L^k)$. Let $\ddbar_k^*\colon \omz^{0,q+1}(M,L^k)\to \omz^{0,q}(M,L^k)$ be the formal adjoint of $\dbar_k$ with respect to $(\,\cdot\, ,\,\cdot\,)_{{k}}$. The Kodaira Laplacian is defined as $\Box_k^q=\ddbar_k^*\ddbar_k+\ddbar_k\ddbar_k^*$ on $\Omega^{0,q}(M,L^k)$. 
We continue to use $\dbar_k$ to denote its maximal extension by
\[
\dbar_{k}: {\rm Dom\,}\dbar_{k}\subset L^2_{(0,q)}(M,L^k)\To L^2_{(0,q+1)}(M,L^k),
\]
with domain ${\rm Dom\,}\dbar_{k}=\big\{u\in L^2_{(0,q)}(M,L^k) \mid  \dbar_{k}u\in L^2_{(0,q+1)}(M,L^k)\big\}$. And let
$$
\dbarstar_{k}:{\rm Dom\,}\dbarstar_{k}\subset L^2_{(0,q+1)}(M,L^k)\to L^2_{(0,q)}(M,L^k)
$$
be the Hilbert space adjoint of $\dbar_{k}$ with respect to $(\,\cdot\, ,\,\cdot\,)_{{k}}$. The Gaffney extension of Kodaira Laplacian is then given by
\begin{equation}\label{e-gue210228yyd}
\begin{split}
&\Box^q_{k}=\dbarstar_{k}\dbar_{k}+\dbar_{k}\dbarstar_{k}:{\rm Dom\,}\Box^q_{k}\subset L^2_{(0,q)}(M,L^k)\To L^2_{(0,q)}(M,L^k),\\
& {\rm Dom\,}\Box^q_{k}=\set{u\in{\rm Dom\,}\dbar_{k}\cap{\rm Dom\,}\dbarstar_{k} \,\big|\, \dbar_{k}u\in{\rm Dom\,}\dbarstar_{k}, \dbarstar_{k} u\in{\rm Dom\,}\dbar_{k}}.
\end{split}
\end{equation}	

For $t>0$, let
$e^{-t\Box^q_{k}}: L^2_{(0,q)}(M,L^k)\To L^2_{(0,q)}(M,L^k)$ denote the heat operator of $\Box_k^q$. 
Let  
\[A_k(t,z,z'):=e^{-t\Box^q_{k}}(z,z')\in\cali{C}^\infty(\mathbb R_+\times M\times M, (T^{*0,q}M\otimes L^k)\boxtimes(T^{*0,q}M\otimes L^k)^*)\]
to be the smooth kernel of $e^{-t\Box^q_{k}}$ with respect to $dv_M$. 
We also write $A_k(t):=e^{-t\Box^q_{k}}$. 	
For $k\in \N$, $A_k(\frac{t}{k}):=e^{-\frac{t}{k}\Box^q_{k}}$ satisfies that
\begin{equation*}%\label{e-gue210531ycda}
\begin{cases}
\Big(\frac{\pa}{\pa t}+\frac{1}{k}\Box^q_{k}\Big)A_k(\frac{t}{k})=0,\\
\lim\limits_{t\To0}A_k(\frac{t}{k})=I\ \ \mbox{on $L^2_{(0,q)}(M,L^k)$}.
\end{cases}
\end{equation*}
From the canonical identification ${\rm End\,}(L^k)=\mathbb C$, we have
\begin{equation}\label{e-gue210303yyd}
A_k(\frac{t}{k},z,z')\in\cali{C}^\infty(\mathbb R_+\times M\times M, T^{*0,q}M\boxtimes(T^{*0,q}M)^*).
\end{equation}

Our main goal of this paper is to study the asymptotic behavior of 
$A_k(\frac{t}{k},z,z)$ as $k\To+\infty$ for $z\in M$. Notice that we have the following unitary identification:
\[
\begin{split}
U: L^2_{(0,q)}(D,L^k)&\longleftrightarrow L^2_{(0,q)}(D),\\
s^k u&\longleftrightarrow e^{-\frac{k\phi}{2}}  u.
\end{split}
\]
For every $\hat u:=s^k u\in\Omega^{0,q}_c(D,L^k)$ with $u\in\Omega^{0,q}_c(D)$, there exists $A_{k,s}(t,z,z')\in\cali{C}^\infty(\mathbb R_+\times D\times D, T^{*0,q}M\boxtimes(T^{*0,q}M)^*)$ such that 
\begin{equation*}
(e^{-t\Box^q_{k}}\hat u)(z)=s^k(z)\, e^{\frac{k\phi(z)}{2}}\int_DA_{k,s}(t,z,z')e^{-\frac{k\phi(z')}{2}} u(z')dv_M(z')\ \ \mbox{on $D$}. 
\end{equation*}
Note that 
\begin{equation}\label{e-gue210303yydII}
A_k(t,z,z')=s^k(z)\,e^{\frac{k\phi(z)}{2}}A_{k,s}(t,z,z')e^{-\frac{k\phi(z')}{2}}(s^k(z'))^*\ \ \mbox{on $D$}. 
\end{equation}
In particular, 
\begin{equation*}%\label{e-gue210303yydIII}
\mbox{$e^{-t\Box^q_{k}}(z,z)=A_k(t,z,z)=A_{k,s}(t,z,z)$ on $D$}.
\end{equation*}

In what follows, we will define a weighted Laplacian $\Box^q_{k\phi}$ and identify $A_{k,s}(t,z,z')e^{\frac{k\phi(z)-k\phi(z')}{2}}$ as the ``distribution kernel" of its heat operator. For $u,v\in \omz_c^{0,q}(D)$, let
$$
(u ,  v)_{k\phi}=(u ,  v)_{k\phi,D}:=\int_{D}\langle\,u ,  v\,\rangle e^{-k\phi}dv_M(z)
$$
denotes the $L^2$-inner product on $\omz_c^{0,q}(D)$ with the weight $k\phi$. Similarly, let $L_{(0,q)}^2(D,k\phi)$ denote the completion of $\omz_c^{0,q}(D)$ with respect to $\|\cdot\|^2_{k\phi}=(\,\cdot\, ,\, \cdot\,)_{k\phi}$, and $\dbar^{*,k\phi}$ the formal adjoint of $\dbar$ with respect to $(\,\cdot\, ,\, \cdot\,)_{k\phi}$. Given point $p\in M$, let $Z_1,\cdots,Z_n$ be an orthonormal frame of $T^{1,0}M$, defined in a neighborhood $D$ of $p$. The dual frame of $Z_1,\cdots,Z_n$ is represented by $w^1,\cdots,w^n$.  Let $(\ov w^j\wedge)^*:=i_{\ov Z_j}$ denotes the contraction on forms. Using these notations, we can locally express:
	
\begin{equation}\label{dbarrho}
\begin{split}
&\dbar=\sum_{j=1}^{n}\ol w^j\wedge \ol Z_j+\sum_{j=1}^{n}(\dbar \ol w^j)\wedge (\ol w^j\wedge)^*,\\
&\dbar^{*,k\phi}=\sum_{j=1}^{n}(\ol w^j\wedge)^*\big(-Z_j+kZ_j\phi+\alpha_j(z)\big)+\sum_{j=1}^{n}\ol w^j\wedge (\dbar \ol w^j\wedge)^*,
\end{split}
\end{equation}
where $\alpha_j(z)$ is a smooth function on $D$, independent of $k$, for every $j=1,\ldots,n$. Now set  
\begin{equation*}\label{bkphi}
\Box_{k\phi}^q:=\dbar^{*,k\phi}\dbar+\dbar\,\dbar^{*,k\phi}: \Dom(\Box_{k\phi}^q)\subset L_{(0,q)}^2(D,k\phi)\To L_{(0,q)}^2(D,k\phi).
\end{equation*}
For $s^ku\in\omz^{0,q}(D,L^k)$ with $u\in \omz^{0,q}(D)$, it is easy to check that
\begin{equation*}\label{bkkphi}
\Box_{k}^q(s^ku)=s^k\Box_{k\phi}^q u.
\end{equation*}
Let
\begin{equation}\label{e-gue210303yydI}
A_{k\phi}(t,z,z'):=e^{\frac{k\phi(z)}{2}}A_{k,s}(t,z,z')e^{\frac{-k\phi(z')}{2}}\ \ \mbox{on $D$}. 
\end{equation}
For $t>0$, define the continuous operator
$A_{k\phi}(t): \mathscr E'(D,T^{*0,q}M)\To\Omega^{0,q}(D)$ by
\begin{equation*}\label{e-gue210305yydI}
A_{k\phi}(t)u=\int A_{k\phi}(t,z,z')u(z')dv_M(z'),\ \ u\in\Omega^{0,q}_c(D).
\end{equation*}
We can verify that 
$$
A'_{k\phi}(t)+\Box^q_{k\phi}A_{k\phi}(t)=0\quad\text{ and }\quad\lim\limits_{t\To0^+}A_{k\phi}(t)=I.
$$	
	 
\section{Heat kernel asymptotics}\label{estimate}  
	
In this section, we will establish the asymptotic expansion for $A_{k\phi}(t,z,z')$ and provide a proof for our main theorem. To achieve this, we will introduce a scaled Laplacian $\Box^q_{(k)}$ of $\Box_{k\phi}^q$ and its corresponding distribution kernel $A_{(k)}(t,z,z')$ in Section \ref{scas}. We then prove such $A_{(k)}(t,z,z')$ is uniformly bounded in $k$ in the $\mathscr{C}^\infty$-topology and has a uniform convergence to the heat kernel of some deformed Laplacians in $\C^n$ (described in Section \ref{mod}).  These uniform properties are essential and can be derived from G\aa{}rding's inequality.

\subsection{Heat kernel on $\C^n$} \label{mod}
\

\vspace{0.18cm}
 
We begin by computing the heat kernel explicitly in the model case $\C^n$ (see Theorem \ref{kernelhn}). Let $z=(z_1,\cdots,z_n)$ be the coordinates of $\C^n$. Then $\big\{\frac{\partial}{\partial z_j}\mid  j=1,\cdots, n\big\}$ forms an orthonormal basis of $T^{1,0}\C^n$ with respect to the standard Hermitian metric $\langle\,\cdot\, , \,\cdot\,\rangle_{\C^n}$, and $\{dz_j\mid  j=1,\cdots, n\}$ is the dual basis.
Let $\dbar$ be the Cauchy-Riemann operator  on $\C^n$, given by
$$
\dbar=\sum_{j=1}^{n}d\ol z_j\wedge\frac{\partial}{\partial \ol z_j}:\omz^{0,q}(\C^n)\to \omz^{0,q+1}(\C^n).
$$
For $\lambda_j\in\R\setminus \{0\}$, $j=1,\cdots, n$, let
\begin{equation*}\label{phi0}
\begin{aligned}  
\phi_0(z)=\sum_{j=1}^{n}\lambda_{j}z_j\ol z_j.
\end{aligned}
\end{equation*}
Let $\dot R_0^L:T^{1,0}\C^n\To T^{1,0}\C^n$ be the linear map such that
\begin{equation}\label{R0}
\langle\,\dot R_0^L U , V\rangle_{\mathbb C^n}=\pa\dbar\phi_{0}(U,\ol V),\quad U,V\in T^{1,0}\C^n.
\end{equation}
Using the notation convention established in Section \ref{sec:prelim}, we define $(\,\cdot \,, \,\cdot\,)_{\phi_0}$ as the weighted inner product on $\omz_c^{0,q}(\C^n)$ with weight $\phi_0$, which satisfies
\begin{equation*}\label{e-gue210506yyd}
\begin{aligned}
(\,u ,  v\,)_{\phi_0}=\int_{\C^n}\langle\,u ,  v\,\rangle_{\C^n}e^{-\phi_0}dv(z), \quad u,v\in\omz_c^{0,q}(\C^n),
\end{aligned}
\end{equation*}
where $dv(z)=i^ndz_1\wedge d\ol z_1\wedge\cdots\wedge dz_{n}\wedge d\ol z_{n}$. Let $L^2_{(0,q)}(\C^n,\phi_0)$ denote the completion of $\omz_c^{0,q}(\C^n)$ with respect to $(\,\cdot \,, \,\cdot\,)_{\phi_0}$, and let $\|\cdot\|_{\phi_0}$ be its corresponding norm. 
The maximal extension of $\dbar$ is still denoted by
$$
\dbar: {\rm Dom\,}\dbar\subset L^2_{(0,q)}(\C^n,\phi_0)\To L^2_{(0,q+1)}(\C^n,\phi_0)
$$
with ${\rm Dom\,}\dbar=\big\{u\in L^2_{(0,q)}(\C^n,\phi_0): \dbar u\in L^2_{(0,q+1)}(\C^n,\phi_0)\big\}$.   
Let 
$$
\dbar^{*,\phi_0}:{\rm Dom\,}\dbar^{*,\phi_0}\subset L^2_{(0,q+1)}(\C^n,\phi_0)\to L^2_{(0,q)}(\C^n,\phi_0)
$$
be the Hilbert space adjoint of $\dbar_0$ with respect to $(\,\cdot \,, \,\cdot\,)_{\phi_0}$. The Gaffney extension of $\Box^q_{\phi_0}$ is then given by
\begin{equation}\label{boxhnn}
\Box^q_{\phi_0}=\dbar^{*,\phi_0}\dbar+\dbar\,\dbar^{*,\phi_0}: {\rm Dom\,}\Box^q_{\phi_0}\subset L^2_{(0,q)}(\C^n,\phi_0)\To L^2_{(0,q)}(\C^n,\phi_0), 
\end{equation}
with ${\rm Dom\,}\Box^q_{\phi_0}=\{ u\in{\rm Dom\,}\dbar\cap{\rm Dom\,}\dbar^{*,\phi_0}:\dbar u\in{\rm Dom\,}\dbar^{*,\phi_0}, \dbar^{*,\phi_0} u\in{\rm Dom\,}\dbar\}$.
Then we can write
\begin{equation}\label{boxhn}
\begin{aligned}
\Box^q_{\phi_0}=&\sum_{j=1}^{n}\Big(-\frac{\pa}{\pa z_j}+\lambda_{j}\bar z_l\Big)\frac{\pa}{\pa \ol z_j}+\sum_{j=1}^{n}\lambda_{j}d\ol z_j\wedge (d\ol z_j\wedge)^* \quad\mbox{on \,$\Omega^{0,q}_c(\C^n)$.}
\end{aligned}
\end{equation}	
	
Let $e^{-t\Box^q_{\phi_0}}$, $t>0$, be the heat operator for $\Box^q_{\phi_0}$, and
$e^{-t\Box^q_{\phi_0}}(z,z')\in\mathscr D'(\mathbb R_+\times \C^n\times \C^n,T^{*0,q}\C^n\boxtimes(T^{*0,q}\C^n)^*)$
be its distribution kernel  with respect to $(\,\cdot \,, \,\cdot\,)_{\phi_0}$. We have
\[
e^{-t\Box^q_{\phi_0}}(z,z')\in\cali{C}^\infty(\mathbb R_+\times \C^n\times \C^n,T^{*0,q}\C^n\boxtimes(T^{*0,q}\C^n)^*).
\]
To compute $e^{-t\Box^q_{\phi_0}}(z,z')$, we shall construct another Laplacian $\Box_0^q$ on $L^2_{(0,q)}(\C^n)$. Consider the following unitary identification: 
\[
\begin{split}
U: L^2_{(0,q)}(\C^n,\phi_0)&\longleftrightarrow L^2_{(0,q)}(\C^n),\quad\\
e^{\frac{\phi_0}{2}} u&\longleftrightarrow  u.
\end{split}
\]
We define the operator
\[
\dbar_0=\sum_{j=1}^nd\ol z_j\wedge\Big(\frac{\pa}{\pa\ol z_j}+\frac{1}{2}\frac{\pa\phi_0}{\pa\ol z_j}\Big) \quad\mbox{on\, $\omz^{0,q}(\C^n)$.}
\]
For $u\in\omz^{0,q}(\C^n)$, it can be verified that that
$$
\dbar_0 u=e^{-\frac{\phi_0}{2}} \dbar (e^{\frac{\phi_0}{2}} u).
$$
Let $\dbarstar_0:  \Omega^{0,q+1}(\mathbb C^n)\To\Omega^{0,q}(\mathbb C^n)$ be the formal adjoint of $\dbar_0$ with respect to $(\,\cdot \,, \,\cdot\,)_{\mathbb C^n}$. Thus we have 
$$
\dbarstar_0=\sum_{j=1}^n(d\ol z_j\wedge)^*\Big(-\frac{\pa}{\pa z_j}+\frac{1}{2}\frac{\pa\phi_0}{\pa z_j}\Big).
$$
Now let 
$$
\Box_0^q:=\dbar_0\dbarstar_0+\dbarstar_0\dbar_0: {\rm Dom\,}\Box_0^q\subset L^2_{(0,q)}(\mathbb C^n)\To L^2_{(0,q)}(\mathbb C^n)
$$
be the self-adjoint extension with 
${\rm Dom\,}\Box_0^q=\big\{u\in L^2_{(0,q)}(\mathbb C^n) : \Box_0^q u\in L^2_{(0,q)}(\mathbb C^n)\big\}$. Then we have
\begin{equation*}%\label{boxeta2}
\begin{aligned}
\Box_0^q=\sum_{j=1}^{n}\Big(-\frac{\pa}{\pa z_j}+\frac{1}{2}\lambda_{j}\ov z_j\Big)\Big(\frac{\pa}{\pa\ol z_j}+\frac{1}{2}\lambda_{j}z_j\Big)
+\sum_{j=1}^{n}\lambda_{j}d\ol z_j\wedge (d\ol z_j\wedge)^* \quad\mbox{on\, $\omz^{0,q}(\C^n)$.}
\end{aligned}
\end{equation*} 
Moreover, 
\begin{equation}\label{box0}
\Box^q_{\phi_0}=e^{\frac{\phi_0}{2}}\Box_0^q e^{-\frac{\phi_0}{2}}.
\end{equation} 
Let $e^{-t\Box_0^q}$ be the heat operator of $\Box_0^q$ and let $e^{-t\Box_0^q}(z,z')\in\cali{C}^\infty(\mathbb R_+\times\mathbb C^n\times\mathbb C^n, T^{*0,q}\mathbb C^n\boxtimes(T^{*0,q}\mathbb C^n)^*)$ be its distribution kernel.
Let 
\[
\begin{split}
\varTheta_0&:=\sum_{j=1}^{n}\lambda_jd\ov z_j\wedge (d\ov z_j\wedge)^*,\\
\mathcal H&:=-\sum_{j=1}^{2n}\bigg(\frac{\partial}{{\sqrt{2}}\partial x_j}+\frac{1}{2}\partial\dbar\phi_0\Big(\sum^{2n}_{k=1}x_k\frac{\partial}{\partial x_k},\frac{\partial}{{\sqrt{2}}\partial x_j}\Big)\bigg)^2-\sum_{j=1}^n \lambda_j,
\end{split}
\] 
where $z_j=x_{2j-1}+ix_{2j}$, $j=1,\cdots,n$. $\mathcal H$ is a second-order partial differential operator and the heat kernel of its self-adjoint extension on $\C^n$ is well-known. Indeed, the Mehler's formula \cite[Appendix E.2]{MM07} and $\eqref{R0}$ gives that
\begin{equation}\label{mehler}
\begin{aligned}
e^{-t\mathcal H}(z,z')=&\frac{1}{(2\pi)^n}\dfrac{\det {\dot R_0^L}}{\det\big(1-e^{-2t{\dot R_0^L}}\big)}
\exp\bigg\{-\frac{1}{2}\Big\langle\,\frac{{\dot R_0^L}/2}{\tanh(t{\dot R_0^L})}z ,  z\,\Big\rangle_{\mathbb C^n}\\
&-\frac{1}{2}\Big\langle\,\frac{{\dot R_0^L}/2}{\tanh(t{\dot R_0^L})}z' ,  z'\,\Big\rangle_{\mathbb C^n}
+\Big\langle\,\frac{{\dot R_0^L}/2}{\sinh(t{\dot R_0^L})}e^{t{\dot R_0^L}}z ,  z'\,\Big\rangle_{\mathbb C^n}\bigg\}, 
\end{aligned}
\end{equation}
where $R_0^L$ is as in \eqref{R0}. Notice that we can rewrite $\Box^q_0$ as
$$
\Box^q_0=\frac{1}{2}\mathcal H+\varTheta_0.
$$
It then follows that
\begin{equation}\label{e-gue210521yyd}
e^{-t\Box^q_0}(z,z')=e^{-t\varTheta_0}e^{-\frac{t}{2}\mathcal H}(z,z').
\end{equation}

For every $t>0$, we define $P(t,z,z')\in\cali{C}^\infty(\mathbb R_+\times \C^n\times \C^n,T^{*0,q}\C^n\boxtimes(T^{*0,q}\C^n)^*)$ by 
\begin{equation}\label{e-gue210521ycdII}
P(t,z,z'):=e^{\frac{\phi_0(z)}{2}} e^{-t\Box^q_0}(z,z')e^{-\frac{\phi_0(z')}{2}},
\end{equation}
Let
\[
P(t): \Omega^{0,q}_c(\C^n)\To\Omega^{0,q}(\C^n)
\]
be the continuous operator given by 
\begin{equation*}\label{e-gue210521ycdIII}
\big(P(t)u\big)(z)=\int P(t,z,z')u(z')dv(z'),\ \ u\in\Omega^{0,q}_c(\C^n).
\end{equation*} 
For every $u,v\in\Omega^{0,q}_c(\C^n)$ and $t>0$, it follows from \eqref{box0} that
\begin{equation}\label{e-gue210521yyda}
\begin{split}
	%&P(t)u\in\Omega^{0,q}(\C^n)\cap L^2_{(0,q)}(\C^n,\phi_0)\\
P(t)u\in{\rm Dom\,}\Box^q_{\phi_0},\quad P'(t)u+\Box^{q}_{\phi_0}\big(P(t)u\big)=0,\quad\mbox{and }
\lim_{t\To0^{+}}P(t)u=u.
\end{split}
\end{equation}
 We now prove that $P(t)=e^{-t\Box^q_{\phi_0}}$.  From \eqref{e-gue210521yyda}, we have
\begin{equation}\label{pe-gue210504ycd}
\begin{split}
&\big(u ,  e^{-t\Box^q_{\phi_0}}v\big)_{\phi_0}-\big(P(t)u ,  v\big)_{\phi_0}=\int^t_0\frac{\partial}{\partial s}\Bigr(\big(\,P(t-s)u ,  e^{-s\Box^q_{\phi_0}}v\,\big)_{\phi_0}\Bigr)ds\\
&\quad=\int^t_0\big(-P'(t-s)u ,  e^{-s\Box^q_{\phi_0}}v\big)_{\phi_0} ds+\int^t_0\Big(P(t-s)u , -\Box^q_{\phi_0}\big(e^{-s\Box^q_{\phi_0}}v\big)\Big)_{\phi_0} ds\\
&\quad=-\int^t_0\big(P'(t-s)u+\Box^q_{\phi_0}P(t-s)u ,  e^{-s\Box^q_{\phi_0}}v\big)_{\phi_0} ds=0.
\end{split}  
\end{equation}
Note that $e^{-t\Box^q_{\phi_0}}$ is self-adjoint, we obtain $P(t)=e^{-t\Box^q_{\phi_0}}$ and hence
\[
e^{-t\Box^{q}_{\phi_0}}(z,z')=e^{\frac{\phi_0(z)}{2}} e^{-t\Box^q_0}(z,z')e^{-\frac{\phi_0(z')}{2}}.
\]
From this observation, \eqref{mehler} and \eqref{e-gue210521yyd}, we conclude that  
\begin{thm}\label{kernelhn}
The heat kernel of the Laplacian $\Box^q_{\phi_0}$ is given by 
$$
e^{-t\Box^{q}_{\phi_0}}(z,z')=e^{\frac{\phi_0(z)}{2}} e^{-t\varTheta_0}e^{-\frac{t}{2}\mathcal H}(z,z')e^{-\frac{\phi_0(z')}{2}}
$$
for $z,z'\in\C^n$.
In particular, at $z=z'=0$, we have
\begin{equation}\label{e-gue210523yydII}
\begin{split}
e^{-t\Box^q_{\phi_0}}(0,0)&=\frac{1}{(2\pi)^n}\dfrac{\det {\dot R_0^L}}{\det\big(1-e^{-t {\dot R_0^L}}\big)}e^{-t\varTheta_0}=\prod_{j=1}^n\frac{\lambda_j\big(1+(e^{-t\lambda_j}-1)d\ov z_j\wedge (d\ov z_j\wedge)^*\big)}{2\pi (1-e^{-t\lambda_j})}.
\end{split}
\end{equation} 
\end{thm}
\begin{rem}\label{aaaa1}
For the case where $\lambda_j = 0$ for all $1 \leq j \leq n$, we have the following result:
\[
e^{-t\Box^q_{\phi_0}}(0,0) = e^{-\frac{t}{2} \big( \sum_{j=1}^{2n} \frac{\partial}{\sqrt{2}\partial x_j} \big)^2}(0,0) = e^{-\frac{t}{2} \Delta}(0,0) = \frac{1}{(2\pi t)^n}.
\]
Here $\big\{ \frac{\partial}{\sqrt{2}\partial x_j} \big\}_{j=1}^{2n}$ is an orthonormal basis of the tangent space $T\mathbb{C}^n$.
\end{rem}

\subsection{Scaling technique} \label{scas}
		\
	
\vspace{0.18cm}

Let's return to the complex manifold case. Consider a point $p \in M$, where $M$ is a complex manifold. Let $\{ Z_j \}_{j=1}^n$ be a local orthonormal frame of $T^{1,0} M$ and $\{ w^j\}_{j=1}^n$ its dual frame, defined on an open neighborhood $D \subset M$ of $p$. Let $s$ be a local holomorphic trivializing section of $L$ on $D$ with local weight $\phi$, i.e., $|s(z)|^2_{h^L}=e^{-\phi(z)}$ for all $z\in D$.

To apply the scaling technique, we need to introduce a suitable local holomorphic coordinate system. For a given $p\in M$, there exists a local holomorphic coordinate $z=(z_1,\cdots,z_n)$ on an open neighborhood $D\subset M$, identified with an open subset in $\C^n$, such that the following properties hold for $k,j=1,\cdots,n$:
\begin{equation}\label{local_co}
\begin{split}
&z(p)=0,\quad
\Big\langle\frac{\pa}{\pa z_j},\frac{\pa}{\pa z_k}\Big\rangle=\delta_{jk}+O(|z|),\\
&Z_j=\frac{\pa}{\pa z_j}+\sum_{s=1}^nr_{j,s}(z)\frac{\pa}{\pa z_s}, \quad r_{j,s}(z)=O(|z|)\in\mathscr C^\infty(D),\\
&\phi(z)=\sum_{j=1}^{n}\lambda_{j}|z_j|^2+O(|z|^3),\quad \lambda_j\in \R.
\end{split}
\end{equation}
Here, the notation $O(r)$ indicates that the term is bounded by $Cr$ for certain constant $C$ as $r\rightarrow 0$.
Note that $z_j=x_{2j-1}+ix_{2j}$, $\frac{\pa}{\pa z_j}=\frac{1}{2}\Big(\frac{\pa}{\pa x_{2j-1}}-i\frac{\pa}{\pa x_{2j}}\Big)$ and $\big\langle\,\frac{\partial}{\partial x_k},\frac{\partial}{\partial x_j}\,\big\rangle=2\delta_{kj}+O(|x|)$. In what follows, we will work with this local coordinate system.

Let $r\in\R$. Define the ball $B_r:=\{z\in\C^{n}, |z|<r\}$, which can be identified as a subset of $M$ via the scaling map:
\begin{equation}\label{eq_scaling}
\begin{split}
F_k: \mathbb C^{n}\To\mathbb C^{n},\quad
F_kz=\frac{z}{\sqrt{k}}.
\end{split}
\end{equation}
Observe that we can choose sufficient large $k$ such that $F_k(B_{\log k})\subset D$. For $z\in B_{\log k}$, let
\begin{equation}
T^{*0,q}_{(k),z}M:=\Bigg\{\,\sumprime_{|J|=q}a_J\ol w^J_{(k)}\,\bigg|\, a_J\in\C\Bigg\},
\end{equation}
where $w^j_{(k)}(z):= w^j(F_k z)$ for $1\le j\le n$, and $\sumprime$ indicating that the sum is taken over strictly increasing multi-indices.
Let $T_{(k)}^{*0,q}M$ be the vector bundle over $B_{\log k}$ with fiber $T^{*0,q}_{(k),z}M$ at $z\in B_{\log k}$. Define $\omz^{0,q}_{(k)}(B_r)$ as the space of smooth sections of $T^{*0,q}_{(k)}M$ over $B_r$, and $\Omega^{0,q}_{(k),c}(B_r)$ as a subspace of $\omz^{0,q}_{(k)}(B_r)$ consisting of elements with compact support in $B_r$.
For $u=\sumprime_{|J|=q}u_J\ol w^J\in \omz^{0,q}(F_k(B_{\log k}))$, the scaled form $u_{(k)}(z)$ is defined by
\begin{equation}
u_{(k)}(z):=\sumprime_{|J|=q}u_J(F_kz)\ol w^J_{(k)}\in \omz^{0,q}_{(k)}(B_{\log k}).
\end{equation}
	
For a fixed $q\in\set{0,1,\ldots,n}$, note that the set
$\set{\overline w^J(z) ,  J=(j_1,\ldots,j_q), 1\leq j_1<\cdots<j_q\leq n}$ 
forms an orthonormal basis for $T^{*0,q}_zM$ at each $z\in D$. Let $\langle\,\cdot\, , \,\cdot\,\rangle_{(k)}$   
be the Hermitian metric on $T^{*0,q}_{(k)}M$ over $B_{\log k}$, such that $\set{\overline w^J_{(k)} ,  J=(j_1,\ldots,j_q), 1\leq j_1<\cdots<j_q\leq n}$ forms an orthonormal frame for every $z\in B_{\log k}$. We write the volume form  as
$dv_M(z)=m(z)d\sigma(z)$ on $D$, where $d\sigma(z):=d x_1\wedge\cdots\wedge dx_{2n}$ and $m(z)\in\cali{C}^\infty(D)$. Let $(\,\cdot \,, \,\cdot\,)_{(k)}=(\,\cdot \,, \,\cdot\,)_{k\phi_{(k)}}$ be the $L^2$-inner product with weight $k\phi_{(k)}$ on $\Omega^{0,q}_{(k),c}(B_{\log k})$, given by 
\[
(\,u ,  v\,)_{(k)}=\int_{B_{\log k}}\langle\,u ,  v\,\rangle_{(k)}e^{-k\phi_{(k)}}m(F_kz)d\sigma(z),\ \ u, v\in \Omega^{0,q}_{(k),c}(B_{\log k}).
\]
The scaled differential operator $\dbar_{(k)}:\omz^{0,q}_{(k)}(B_{\log k})\to \omz^{0,q+1}_{(k)}(B_{\log k})$ is defined as
\begin{equation}\label{dbarrhok}
\begin{aligned}
\dbar_{(k)}&=\sum_{j=1}^{n}\ol w^j_{(k)}\wedge \ol Z_{j,(k)}+\sum_{j=1}^{n}\frac{1}{\sqrt k}(\dbar\ol w^j)_{(k)}\wedge \big(\ol w^j_{(k)}\wedge\big)^*.
\end{aligned}
\end{equation}
where
$$
\ol Z_{j,(k)}:=\frac{\pa}{\pa \ol z_j}+\sum_{s=1}^n \ov r _{j,s}(F_kz)\frac{\pa}{\pa \ol z_s},
$$
Compare \eqref{dbarrhok} with \eqref{dbarrho}, we observe that
\begin{equation}\label{k1}
\dbar_{(k)}u_{(k)}=\frac{1}{\sqrt k}(\dbar u)_{(k)},\quad\mbox{ for all }u\in\omz^{0,q}(F_k(B_{\log k})).
\end{equation} 
Let $\dbar_{(k)}^*$ be the formal adjoint of $\dbar_{(k)}$ with respect to $(\,\cdot\, , \,\cdot\,)_{(k)}$, then
\begin{equation}\label{dbarstarrhok}
\begin{aligned}
\dbar_{(k)}^*=&\sum_{j=1}^{n}\big(\ol w^j_{(k)}\wedge\big)^*\Big(-Z_{j,(k)}+\sqrt k(Z_j\phi)_{(k)}+\frac{1}{\sqrt k}(\alpha_j)_{(k)}\Big)\\
&+\frac{1}{\sqrt k}\ol w^j_{(k)}\wedge \big((\dbar \ol w^j)_{(k)}\wedge\big)^*, 
\end{aligned}
\end{equation}
where $\alpha_j=\alpha_j(z)$ is a smooth function as in \eqref{dbarrho}. Moreover, we have
\begin{equation}\label{k2}
\dbar^{*}_{(k)}u_{(k)}=\frac{1}{\sqrt k}\big(\dbar^{*,k\phi}u\big)_{(k)},\quad\mbox{ for all } u\in\Omega^{0,q+1}(F_k(B_{\log k})).
\end{equation}
We now define the scaled Kodaira Laplacian as
\begin{equation*}\label{e-gue210614yydI}
\Box^q_{(k)}=\dbar_{(k)}^*\dbar_{(k)}+\dbar_{(k)}\dbar_{(k)}^*
\end{equation*}
on $\Omega_{(k)}^{0,q}(B_{\log k})$.  Then \eqref{k1} and \eqref{k2} yield that
\begin{equation}\label{k(k)}
\Box_{(k)}^qu_{(k)}=\frac{1}{k}(\Box_{k\phi}^qu)_{(k)},\quad\mbox{ for all } u\in\Omega^{0,q}(F_k(B_{\log k})). 
\end{equation} 
	
For $z,z'\in B_{\log k}$, let
\begin{equation}\label{e-gue210325yyd}
\begin{split}
A_{(k)}(t,z,z'):=k^{-n}A_{k\phi}(\frac{t}{k},F_kz,F_kz'),
\end{split}
\end{equation}
where $F_kz=\frac{z}{\sqrt{k}} $, $F_kz'=\frac{z'}{\sqrt{k}}$. 
Let  $A_{(k)}(t): \Omega^{0,q}_{(k),c}(B_{\log k})\To \Omega^{0,q}_{(k)}(B_{\log k})$
be the continuous operator given by 
\begin{equation*}\label{e-gue210325yydI}
\big(A_{(k)}(t)u\big)(z)=\int A_{(k)}(t,z,z')u(z')m(F_kz')d\sigma(z'),\ \ u\in \Omega^{0,q}_{(k),c}(B_{\log k}). 
\end{equation*} 
From \eqref{dbarrhok}, we have
\begin{equation}\label{staru3} 
\ol Z_{j,(k)}=\frac{\pa}{\pa \ol z_j}+\epsilon_kU_{j,k},
\end{equation}
where $\epsilon_k\to 0$ as $k\to+\infty$, and $U_{j,k}$ is a first-order differential operator with coefficients whose  derivatives are uniformly bounded in $k$ on $B_{\log k}$.
Specifically,  $U_{j,k}=\sum_{s=1}^{n}f_{j,s}(k,z)\frac{\pa}{\pa \ov z_s}$ satisfying, for each $1\leq s\leq n$, $l\in \N$, there exits a constant $C_l$ independent of $k$ and $z$ such that 
\begin{equation}\label{ubink}
\sup_{|\alpha|=l,~z\in B_{\log k}}\left |\partial_z^{\alpha} f_{j,s}(k,z) \right|<C_l.
\end{equation}
On the other hand, from the last equality in \eqref{local_co}, we see that in \eqref{dbarstarrhok},
\begin{equation*}\label{staru4}
\begin{aligned}
-Z_{j,(k)}+\sqrt k(Z_j\phi)_{(k)}+\frac{1}{\sqrt k}(\alpha_j)_{(k)}
=-\frac{\pa}{\pa z_j}+\lambda_{j}\bar z_j+\delta_kV_{j,k},
\end{aligned}  
\end{equation*}
where $\delta_k\to 0$ as $k\to+\infty$, and $V_{j,k}$ is a first-order differential operator with coefficients whose  derivatives are uniformly bounded in $k$ on $B_{\log k}$, as in \eqref{ubink}. Combining this with \eqref{dbarrhok}, \eqref{dbarstarrhok}, and \eqref{staru3}, we obtain the following proposition:
	
\begin{prop}\label{3.3}
We have that
\begin{equation}\label{boxrhok}
\begin{aligned}
\Box_{(k)}^{q}=&\sum_{j=1}^{n}\Big(-\frac{\pa}{\pa z_j}+\lambda_{j}\bar z_j\Big)\frac{\pa}{\pa \ol z_j}+\sum_{j=1}^{n}\lambda_{j}\,\ol w^j_{(k)}\wedge \big( \ol w^j_{(k)}\wedge\big)^*+\varepsilon_kP_k,
\end{aligned}
\end{equation}
where $\varepsilon_k$ is a sequence tending to zero as $k\to+\infty$ and $P_k$ is a second-order differential operator with coefficients  uniformly bounded in $k$ on $B_{\log k}$.
\end{prop}

Let $L^2_{(k)}(B_{\log k},T_{(k)}^{*0,q}M)$ be the completion of $\Omega^{0,q}_{c,(k)}(B_{\log k})$ with respect to $(\,\cdot \,, \,\cdot\,)_{(k)}$. The maximal extension of $\dbar_{(k)}$ is still denoted by
\[\dbar_{(k)}: {\rm Dom\,}\dbar_{(k)}\subset L^2_{(k)}(B_{\log k},T_{(k)}^{*0,q}M)\To L^2_{(k)}(B_{\log k},T_{(k)}^{*0,q+1}M),\]
with $ {\rm Dom\,}\dbar_{(k)}=\big\{u\in L^2_{(k)}(B_{\log k},T_{(k)}^{*0,q}M) \mid \dbar_{(k)}u\in L^2_{(k)}(B_{\log k},T_{(k)}^{*0,q+1}M)\big\}$. 
Let $\dbar^{*}_{(k)}$ be the Hilbert space adjoint of $\dbar_{(k)}$ with respect to $(\,\cdot \,, \,\cdot\,)_{(k)}$. The Gaffney extension of $\Box^q_{(k)}$ is then given by 
\[
\begin{split}
\Box^q_{(k)}=\dbar^{*}_{(k)}\dbar_{(k)}+\dbar_{(k)}\dbar^{*}_{(k)}
\mid {\rm Dom\,}(\Box^q_{(k)})\subset L^2_{(k)}(B_{\log k},T_{(k)}^{*0,q}M)\To L^2_{(k)}(B_{\log k},T_{(k)}^{*0,q}M)
\end{split}
\] 
with ${\rm Dom\,}(\Box^q_{(k)})=\{ 
u\in{\rm Dom\,}\dbar_{(k)}\cap{\rm Dom\,}\dbar^{*}_{(k)}:
\dbar_{(k)}u\in{\rm Dom\,}\dbar^{*}_{(k)}, 
\dbar^{*}_{(k)}u\in{\rm Dom\,}\dbar_{(k)}\}$.

For $s\in\mathbb N$ and $B\subset B_{\log k}$, define 
$W^s_{(k)}(B, T_{(k)}^{*0,q}M)$
as the $L^2$-Sobolev space of order $s$ on the sections of $T_{(k)}^{*0,q}M$ over $B$ with respect to $(\,\cdot \,, \,\cdot\,)_{(k)}$, with Sobolev norm 
\begin{equation}
\begin{aligned}
\|u\|^2_{(k),s,B}=\sum_{|\alpha|\le s}\sumprime_{|J|=q}\int_{B}|\pa_{z}^\alpha u_J|^2e^{-k\phi_{(k)}}m(F_kz)d\sigma(z),
\end{aligned}
\end{equation}
for $u=\sumprime_{|J|=q}u_J\ol w^j_{(k)}\in W^s_{(k)}(B, T_{(k)}^{*0,q}M)$. The $L^2$-norm is denoted by $\|\cdot\|_{(k),B}:=\|\cdot\|_{(k),0,B}$.
We define the spaces
\[
\begin{split}
&W^s_{(k),{\rm c\,}}(B, T_{(k)}^{*0,q}M):=\set{u\in W^s_{(k)}(B, T_{(k)}^{*0,q}M) : {\rm supp\,}u\Subset B},\\
&W^s_{(k),{\rm loc}}(B, T_{(k)}^{*0,q}M):=\set{u\in\mathscr D'(B,T_{(k)}^{*0,q}M) : \mbox{$\chi u\in W^s_{(k)}(B, T_{(k)}^{*0,q}M)$, for every $\chi\in\cali{C}^\infty_c(B)$}}.
\end{split}
\]
For $s\in\mathbb Z$, $s<0$, the spaces $W^{s}_{(k),{\rm c\,}}(B, T_{(k)}^{*0,q}M)$ and $W^s_{(k),{\rm loc}}(B, T_{(k)}^{*0,q}M)$ are defined as the dual spaces of $W^{-s}_{(k),{\rm loc}}(B, T_{(k)}^{*0,q}M)$ and $W^{-s}_{(k),{\rm c}}(B, T_{(k)}^{*0,q}M)$, respectively, with respect to $(\,\cdot \,, \,\cdot\,)_{(k)}$. 
\vspace{0.15cm}
	 
The following basic elliptic estimate can be derived from the ellipticity of $\Box_{(k)}^q$ on $B_{2r}$ 
and Proposition \ref{3.3}:
	
\begin{prop}\label{p-gue210511yyd}  
Let $s\in\mathbb N$. For $r>0$ with $B_{2r}\subset B_{\log k}$, there exists a constant $C_{r,s}>0$, independent of $k$, such that 
\begin{equation}\label{subellipticest}
\begin{aligned}
\|u\|^2_{(k),s+2,B_r}\le C_{r,s}\Big(\|u\|^2_{(k),B_{2r}}+\|\Box_{(k)}^qu\|^2_{(k),s,B_{2r}}\Big)
\end{aligned}
\end{equation} 
 for all $u\in \omz_{(k)}^{0,q}(B_{\log k})$.
\end{prop}
Using \eqref{subellipticest}, we can show that $A_{(k)}(t,z,z')$ is uniformly bounded in $k$ in the $\cali{C}^\infty$ topology over $B_{\log k} \times B_{\log k}$ for a fixed $t > 0$.  To establish this, note that the Kodaira Laplacian $\Box^q_k$, defined in \eqref{e-gue210228yyd}, is self-adjoint. By the spectral theorem \cite{D95}, the space $L^2_{(0,q)}(M, L^k)$ can be identified with $L^2(\mathbb{S} \times \mathbb{N}, d\mu)$, and $\Box^q_k$ corresponds to the multiplication operator $M_s$. Here, $\mathbb{S}$ represents the spectrum of $\Box^q_k$, and $\mu$ is a finite regular Borel measure on $\mathbb{S} \times \mathbb{N}$.  In particular, the heat operator $e^{-t\Box^q_k}$ acts as follows:
\[
\begin{aligned}
e^{-t\Box^q_k}: L^2(\mathbb{S} \times \mathbb{N}, d\mu) &\to L^2(\mathbb{S} \times \mathbb{N}, d\mu), \\
\varphi(s,n) &\mapsto e^{-ts}\varphi(s,n).
\end{aligned}
\]  
	
\begin{lem}\label{l-gue210331yyd}
Fix $t>0$. For every $\ell\in\mathbb N$ and $u\in\Omega^{0,q}(M,L^k)\cap \Dom(\Box^q_k)$, we have 
\begin{equation}\label{e-gue210331yyd}
\norm{e^{-\frac{t}{k} \Box^q_{k}}( \Box^q_{k})^\ell u}_{k}=\norm{( \Box^q_{k})^\ell e^{-\frac{t}{k} \Box^q_{k}}u}_{k}
\leq(1+\frac{k^\ell}{t^\ell})C_\ell\norm{u}_{k}, 
\end{equation}
where $C_\ell>0$ is a constant independent of $k$ and $t$. 
\end{lem} 
	
\begin{proof}
Let $\ell\in\mathbb N$. For $u\in\Omega^{0,q} (M,L^k)\cap \Dom(\Box^q_k)$, we identify it with $\varphi(s,n)\in L^2(\mathbb S\times\mathbb N,d\mu)$, then
\[
\begin{split}
&\|e^{-\frac{t}{k} \Box^q_{k}}( \Box^q_{k})^\ell u\|^2_{k}=\|( \Box^q_{k})^\ell e^{-\frac{t}{k} \Box^q_{k}}u\|^2_{k}\\
&=\int\abs{e^{-\frac{t}{k}s}s^\ell \varphi(s,n)}^2d\mu=\int_{\set{0\leq s<1}}\abs{e^{-\frac{t}{k}s}s^\ell \varphi(s,n)}^2d\mu+
\int_{\set{s\geq1}}\abs{e^{-\frac{t}{k}s}s^\ell \varphi(s,n)}^2d\mu\\
&\leq\int_{\set{0\leq s<1}}\abs{\varphi(s,n)}^2d\mu+
\hat C_\ell \int_{\set{s\geq1}}\abs{\frac{k^\ell }{t^\ell s^\ell }s^\ell \varphi(s,n)}^2d\mu\\
&\leq(1+\hat C_\ell \frac{k^{2\ell }}{t^{2\ell }})\int\abs{\varphi(s,n)}^2d\mu,
\end{split}
\]
where $\hat C_\ell >0$ is a constant independent of $k$ and $t$. The lemma follows. 
\end{proof}
We can now prove 
\begin{prop}\label{uniformly}
Let $I\subset\mathbb R_+$ be a compact set and $t\in I$.
For $\ell\in\mathbb N$ and $r>0$ with $B_{r}\subset B_{\log k}$, there exists a constant $C_{\ell,r}>0$, independent of $k$ and $t$, such that 
\begin{equation}\label{uniest}
\|A_{(k)}(t,z,z')\|_{\cali{C}^\ell(I\times B_r\times B_r, T_{(k)}^{*0,q}M\boxtimes(T_{(k)}^{*0,q}M)^*)}\leq C_{\ell,r}.
\end{equation}
\end{prop}
	
\begin{proof}		
Let $\chi\in\cali{C}^\infty_c(B_r,\mathbb R)$ be a non-negative cut-off function.
We claim that $A_{(k)}(t)$ can be continuously extended to $W^{-\ell}_{(k),{\rm c}}(B_r,T_{(k)}^{*0,q}M)$, and 
\begin{equation}\label{e-gue210412yyda}
A_{(k)}(t)\chi: W^{-\ell}_{(k),{\rm c}}(B_r,T_{(k)}^{*0,q}M)\To W^{s}_{(k),{\rm loc}}(B_r,T_{(k)}^{*0,q}M)
\end{equation}
is continuous for all $s\in\mathbb Z$, with the continuity being uniform in $t\in I$.
		
Consider $(\Box^q_{(k)})^\ell A_{(k)}(t)$. Let $u\in \Omega^{0,q}_{(k),c}(B_r)$ and $v\in\Omega^{0,q}_c(F_k(B_r))$ such that $u=v_{(k)}$ on $B_r$. On $D$, we identify $ \Box^q_{k}$ with $\Box^q_{k\phi}$ and sections of $L^k$ with functions. From \eqref{k(k)} and \eqref{e-gue210325yyd}, we have
\begin{equation*}\label{e-gue210410yyd}
\begin{split}
&(\Box^q_{(k)})^{\ell}A_{(k)}(t)u=(\Box^q_{(k)})^{\ell}\Bigr(A_{k\phi}(\frac{t}{k})v\Bigr)_{(k)}=k^{-\ell}\Bigr((\Box^q_{k\phi})^\ell A_{k\phi}(\frac{t}{k})v\Bigr)_{(k)} \,
\, \mbox{ on $B_r$}.
\end{split}
\end{equation*} 
From this and Lemma~\ref{l-gue210331yyd}, we obtain
\begin{equation}\label{e-gue210412yyd}
\begin{split}
&\|(\Box^q_{(k)})^\ell A_{(k)}(t)u\|^2_{(k),B_r}=\Big\|k^{-\ell}\Bigr((\Box^q_{k\phi})^\ell A_{k\phi}(\frac{t}{k})v\Bigr)_{(k)}\Big\|^2_{(k),B_r}\\
&=k^{-2\ell+n}\Big\|(\Box^q_{k\phi})^\ell A_{k\phi}(\frac{t}{k})v\Big\|^2_{k\phi,F_k(B_r)}
= k^{-2\ell+n}\big\|( \Box^q_{k})^\ell e^{-\frac{t}{k} \Box^q_{k}}(s^kv)\big\|^2_{k}\\
&\leq\frac{k^{n}}{t^{2\ell}}C_\ell\|s^kv\|^2_{k}=\frac{k^{n}}{t^{2\ell}}C_\ell\|v\|^2_{k\phi,F_k(B_r)}
=\frac{C_\ell}{t^{2\ell}}\|u\|^2_{(k),B_r}\leq\hat C_\ell\|u\|^2_{(k),B_r}
\end{split}
\end{equation}
for every $t\in I$, where $C_\ell>0$, $\hat C_\ell>0$ are constants independent of $k$ an $t$. From Proposition~\ref{p-gue210511yyd} and \eqref{e-gue210412yyd}, we conclude that 
\begin{equation}\label{e-gue210412yydI}
A_{(k)}(t): L^2_{(k),{\rm c}}(B_r,T_{(k)}^{*0,q}M)\To W^{2s}_{(k),{\rm loc}}(B_r,T_{(k)}^{*0,q}M)
\end{equation}
is continuous for every $s\in\mathbb Z$, with the continuity being uniform in $t\in I$. 
The claim then follows by applying the argument used to prove \eqref{e-gue210412yydI} to the operator $A_{(k)}(t)(\Box^q_{(k)})^{\ell_1}$ for any $\ell_1\in \mathbb N$, taking the adjoint and combining it with Proposition \ref{p-gue210511yyd} and the fact that 
\[
(\Box^q_{(k)})^\ell: L^2_{(k),{\rm c}}(B_r,T_{(k)}^{*0,q}M)\To W^{-2\ell}_{(k),{\rm c}}(B_r,T_{(k)}^{*0,q}M).
\]   
is continuous, see \cite[(3.44)-(3.49)]{HZ23} for details. From \eqref{e-gue210412yyda} and Sobolev embedding theorem, we conclude the proposition.		 
\end{proof}

\begin{thm}\label{t-gue210503yyd}
Let $I\subset\mathbb R_+$ be a bounded interval, and let $r>0$ satisfy $B_{r}\subset B_{\log k}$. We have 
\[
\lim_{k\To+\infty}A_{(k)}(t,z,z')=e^{-t\Box^q_{\phi_0}}(z,z')
\]
in $\cali{C}^\infty(I\times B_r\times B_r,T^{*0,q}\C^n\boxtimes(T^{*0,q}\C^n)^*)$ topology.
\end{thm} 
	
\begin{proof}
From Proposition~\ref{uniformly} and the Cantor diagonal argument, we can extract a subsequence $\{k_1<k_2<\cdots\}$ of $\mathbb N$, with $\lim\limits_{j\To+\infty}k_j=+\infty$, such that 
\[
\lim_{j\To+\infty}A_{(k_j)}(t,z,z')=Q(t,z,z')
\]
locally uniformly on $\mathbb R_+\times \C^n\times \C^n$ in $\cali{C}^\infty$ topology, where 
$Q(t,z,z')\in\cali{C}^\infty(\mathbb R_+\times \C^n\times \C^n,T^{*0,q}\C^n\boxtimes(T^{*0,q}\C^n)^*)$. Consider the continuous operator
$Q(t): \Omega^{0,q}_c(\C^n)\To\Omega^{0,q}(\C^n)$ given by 
\[
(Q(t)u)(z)=\int Q(t,z,z')u(z')dv(z'),\ \ u\in\Omega^{0,q}_c(\C^n).
\]
We claim that for every $u\in\Omega^{0,q}_c(\C^n)$ and $t>0$, we have
\begin{equation}\label{e-gue210503ycd}
\begin{split}
Q(t)u\in{\rm Dom\,}\Box^q_{\phi_0},\quad Q'(t)u+\Box^q_{\phi_0}Q(t)u=0, \quad \mbox{and } \lim_{t\To0^{+}} Q(t)u=u.
\end{split}
\end{equation}
Note that
\begin{equation}\label{e-gue210504yyd}
A'_{(k)}(t)+\Box^q_{(k)}A_{(k)}(t)=0\ \ \mbox{on $B_{\log k}$}.
\end{equation}
Using \eqref{boxrhok} and passing to the limit $k\To+\infty$ in \eqref{e-gue210504yyd}, we obtain $Q'(t)u+\Box^q_{\phi_0}Q(t)u=0$. Write $u\in\Omega^{0,q}_c(\C^n)$ as $u=\sum'_{\abs{J}=q}u_J(z)d\overline z_J$ with $u_J\in\cali{C}^\infty_c(\C^n)$. We set
\begin{equation}\label{e-gue210504yydbc}
u_k:=\sideset{}{'}\sum_{\abs{J}=q}u_J(z)\overline w^J_{(k)}, \quad k=1,2,\cdots.
\end{equation}
Let $r>0$ be fixed. For every $\ell\in\mathbb N$ and $t>0$, we can verify that
\begin{equation}\label{e-gue210504yydI}
\|(\Box^q_{\phi_0})^\ell Q(t)u\|_{\phi_0,B_r}=\lim_{j\To+\infty}\|(\Box^q_{(k_j)})^\ell A_{(k_j)}(t)u_{k_j}\|_{(k_j),B_r}.
\end{equation}
From \eqref{e-gue210412yyd} and \eqref{e-gue210504yydI}, it follows that there is a constant $C_\ell>0$, independent of $t$ and $r$, such that 
\begin{equation}\label{e-gue210504yydII}
\|(\Box^q_{\phi_0})^\ell Q(t)u\|_{\phi_0,B_r}\leq\frac{C_\ell}{t^\ell}\|u\|_{\phi_0,B_r}. 
\end{equation}
Choose $r\gg1$ so that ${\rm supp\,}u\subset B_r$. Then \eqref{e-gue210504yydII} implies
$\|(\Box^q_{\phi_0})^\ell Q(t)u\|_{\phi_0,B_r}\leq\frac{C_\ell}{t^\ell}\|u\|_{\phi_0}$ for every $r\gg1$. Letting $r\To+\infty$, we get
\begin{equation*}\label{e-gue210504yydIII}
\|(\Box^q_{\phi_0})^\ell Q(t)u\|_{\phi_0}\leq\frac{C_\ell}{t^\ell}\|u\|_{\phi_0},
\end{equation*}
which gives  $Q(t)u\in{\rm Dom\,}\Box^q_{\phi_0}$.
To complete the claim, it remains to prove the last statement in \eqref{e-gue210503ycd}. Let $u\in\Omega^{0,q}_c(\C^n)$ and let $u_k\in \Omega^{0,q}_{(k),c}(B_{\log k})$ be as defined in \eqref{e-gue210504yydbc}. For every $t>0$, we have
\begin{equation}\label{e-gue210504yydb}
A_{(k)}(t)u_k-u_k=\int^t_0A'_{(k)}(s)u_kds.%=-\int^t_0\Box^q_{(k)}(A_{(k)}(s)u_k)ds.
\end{equation}
From \eqref{e-gue210412yyd}, we know there exist constants $C>0$ and $\hat C>0$, independent of $t$ and $k$, such that for sufficient large $r$,
\begin{equation*}\label{e-gue210504yydc}
\|A'_{(k)}(s)u_k\|_{(k),B_r}\leq C\|u_k\|_{(k),B_r}\leq \hat C.
\end{equation*}
Using the Lebesgue-dominated convergence Theorem, we obtain 
\[
\begin{split}
&Q(t)u-u=\lim_{j\To+\infty}\int^t_0A'_{(k_j)}(s)u_{k_j}ds=\int^t_0\lim_{j\To+\infty}A'_{(k_j)}(s)u_{k_j}ds\\
&=\int^t_0Q'(s)uds=Q(t)u-\lim_{t\To0+}Q(t)u.
\end{split}
\]
Thus, the claim is concluded.  
		
%We are left with proving that $Q(t)=e^{-t\Box^q_{\phi_0}}$. To do so, we can simply repeat the procedure for operator $Q(t)$ as in \eqref{pe-gue210504ycd}. Based on the preceding discussion, it is clear that for any subsequence of $A_{k}(t)$, we can always find a further sub-subsequence that converges locally uniformly to the same $e^{-t\Box^q_{\phi_0}}$. As a result, we can conclude that $\lim\limits_{k\to+\infty}A_{(k)}(t,z, z')=e^{-t\Box^q_{\phi_0}}( z, z')$ exhibits local uniform convergence on $\mathbb R_+\times M\times M$ in the $\cali{C}^\infty$ topology.

To complete the proof, it remains to establish that $Q(t)=e^{-t\Box^q_{\phi_0}}$. This can be achieved by repeating for $Q(t)$ the same procedure used in \eqref{pe-gue210504ycd}. From the preceding discussion, it follows that for any subsequence of $\{A_{(k)}(t)\}$, there exists a further sub-subsequence that converges locally uniformly to $e^{-t\Box^q_{\phi_0}}$. Consequently, we conclude that $\lim_{k \to +\infty} A_{(k)}(t, z, z') = e^{-t\Box^q_{\phi_0}}(z, z')$ with local uniform convergence on $\mathbb{R}_+ \times M \times M$ in the $\cali{C}^\infty$-topology.
\end{proof}
	
\subsection{Heat kernel asymptotics on $M$} 
\
\vspace{0.18cm}  
	    
We are now able to prove our main results:
	
\noindent\textbf{Proof of Theorem \ref{heatkernelthm}}:	 Inequality \eqref{eq_ub} follows directly from Proposition~\ref{uniformly}. To show \eqref{eq_hka}, let $ p \in M $ and $ 0 \leq q \leq n $, we choose local coordinates $ (z_1, \dots, z_n) $ such that $ z_j(p) = 0 $ for $ j = 1, \dots, n $, and satisfy \eqref{local_co} in a neighborhood of $ p$.
Using Theorems \ref{kernelhn} and \ref{t-gue210503yyd}, along with equations \eqref{e-gue210303yydI} and \eqref{e-gue210325yyd}, we derive
\begin{equation}\label{aaa2}
 \lim_{k\to+\infty}k^{-n}e^{-\frac{t}{k}\Box^q_{k}}(0,0)= e^{-t \Box^q_{\phi_0}}(0,0) = \frac{\det(\dot{R^L}/2\pi) \exp(t\varTheta)}{\det(1 -\exp(-t \dot{R^L}))}(p).   
\end{equation}
Applying this procedure to each point $z \in M$, replacing $0$ with $z$, we establish \eqref{eq_hka}.

 \qed

\begin{rem}
    By applying Theorem \ref{heatkernelthm} to a relatively compact domain $U \Subset M$ and using the dominated convergence theorem, we obtain that for $0 \leq q \leq n$,
    \begin{equation}\label{eq-trlocal} 
    \begin{split}
        &\lim_{k \to +\infty} \int_U \frac{1}{k^n} \Tr_q \big[\exp\big(-\frac{t}{k} \Box_k^q\big)(z,z)\big] dv_M = \int_U \lim_{k \to +\infty} \frac{1}{k^n} \Tr_q \big[\exp\big(-\frac{t}{k} \Box_k^q\big)(z,z)\big] dv_M \\
        &= \int_U \frac{\det(\dot{R^L}/2\pi) \Tr_q \exp(t \varTheta)}{\det(1 -\exp(-t \dot{R^L}))}(z) dv_M \longrightarrow \int_U 1_{U(q)} (-1)^q \frac{c_1(L, h^L)^n}{n!}, \quad t \to \infty.
    \end{split}
    \end{equation}
    %Under additional assumptions on the integration of the heat kernel over $M\setminus U$, new Morse inequalities can be derived by following the approach in \cite[Proposition 3.2.7, Theorem 3.2.13]{MM07}.
\end{rem}

\subsection{Heat kernel asymptotics at degenerate point} \label{sec_rl0}
\
\vspace{0.18cm} 

As noted in Remark \ref{aaaa1}, if an eigenvalue of $\dot{R^L}$ at $p\in M$ is zero in \eqref{aaa2}, then its contribution to the term $\det(\dot{R^L}/2\pi)/\det(1 -\exp(-t \dot{R^L}))$ is $1/(2\pi t)$. Furthermore, if $ \lambda_j(p) = 0 $ for all $ 1 \leq j \leq n $ in the local weight expression $ \phi(z) = \sum_{j=1}^n \lambda_j |z_j|^2 + O(|z|^3) $, we have:
\[
\lim_{k\to+\infty} k^{-n} e^{-\frac{t}{k} \Box^q_{k}}(p,p) = \frac{1}{(2\pi t)^n}.
\]
This raises a natural question: if $\phi$ has a vanishing order of $\sigma$ with $\sigma > 3$ near a degenerate point $p$, can we obtain further insights into the asymptotic expansion of the heat kernel? Specifically, consider the case where $\phi$ can be expressed in a neighborhood of $p = 0$ as:
\[
\phi(z) = \sum_{j=1}^n \lambda_j |z_j|^\sigma + O(|z|^{\sigma+1}), \quad \sigma \geq 3,
\]
with $\lambda_j\neq 0$ for some $j$. To explore this setting, we introduce a new scaling map: 
\begin{equation}\label{aaaa3}
 F_k: \mathbb{C}^{n} \to \mathbb{C}^{n}, \quad F_k z = \frac{z}{k^{1/\sigma}}.   
\end{equation}
Following the procedure outlined in Sec.\,\ref{scas}, we scale the relevant objects. Starting with
\[
\overline{\omega}^j_{(k)}(z) := \overline{\omega}^j(F_k z), \quad u_{(k)}(z) := \sumprime_{|J|=q} u_J(F_k z) \overline{\omega}^J_{(k)} \in \Omega^{0,q}_{(k)}(B_{\log k}),
\]
and
\[
\overline{Z}_{j,(k)} := \frac{\partial}{\partial \overline{z}_j} + \sum_{s=1}^n \overline{r}_{j,s}(F_k z) \frac{\partial}{\partial \overline{z}_s}.
\]
The scaled $ \dbar $-operator becomes
\[
\dbar_{(k)} := \sum_{j=1}^{n} \overline{\omega}^j_{(k)} \wedge \overline{Z}_{j,(k)} + \sum_{j=1}^{n} \frac{1}{k^{1/\sigma}} (\dbar \overline{\omega}^j)_{(k)} \wedge \big(\overline{\omega}^j_{(k)} \wedge \big)^*,
\]
and its adjoint is
$$
\dbar_{(k)}^*:=\sum_{j=1}^{n}\big(\ol w^j_{(k)}\wedge\big)^*\big(-Z_{j,(k)}+k^{1-1/\sigma}(Z_j\phi)_{(k)}+k^{-1/\sigma}(\alpha_j)_{(k)}\big)+\sum_{j=1}^{n}\frac{1}{k^{1/\sigma}}\ol w^j_{(k)}\wedge \big((\dbar \ol w^j)_{(k)}\wedge\big)^*. 
$$
The corresponding scaled Kodaira Laplacian $ \Box^q_{(k)} $ is similarly defined, and importantly, we check that
\[
\Box_{(k)}^q u_{(k)} = k^{-2/\sigma} (\Box_{k\phi}^q u)_{(k)}, \quad \text{for all } u \in \Omega^{0,q}(F_k(B_{\log k})).
\]
Define the scaled heat kernel as
$$
A_{(k)}(t,z,z'):=k^{-2n/\sigma}A_{k\phi}(k^{-2/\sigma}t,F_kz,F_kz')\quad \text{for } z,z'\in B_{\log k}.
$$
Let $ \Box^q_{\phi_0} $ be the Laplacian with $ \phi_0 = \sum_{j=1}^n\lambda_j|z_j|^\sigma $ as defined in \eqref{boxhnn}, and let $ e^{-t \Box_{\phi_0}^q}(z,z') $ denote the heat kernel of $ \Box_{\phi_0}^q $ in $ \mathbb{C}^n $. With the same approach as in Sec.\,\ref{scas}, we obtain the following asymptotic result:
\[
\lim_{k \to \infty} \frac{1}{k^{2n/\sigma}} e^{-\frac{t}{k^{2/\sigma}} \Box_k^q}(p, p)=\lim_{k \to \infty}A_{(k)}(t,0,0) = e^{-t \Box_{\phi_0}^q}(0,0).
\]
This result provides new insights into the heat kernel asymptotics when the local weight has a higher vanishing order. Calculating the heat kernel for $\Box_{\phi_0}^q$ remains challenging. However, as the time parameter $t \to 0$, the heat kernel converges to the Bergman kernel $B_{\phi_0}^{q}$, which is computable in certain cases (such as \cite{MS23}).  Our method offers a simple yet effective way to study the heat kernel, even at degenerate points.

\section{Morse inequalities on coverings}\label{sec-morse}

\subsection{Symmetric tensor power of vector bundles}   
		\
	
\vspace{0.18cm}
The estimates of heat kernel \eqref{eq_ub}-\eqref{eq_hka} are local, thus they can be utilized on possibly non-compact manifolds, for example, covering manifolds with discrete group action. Note that the holomorphic Morse inequalities on coverings had been obtained by Chiose-Marinescu-Todor, which were motivated by Atiyah's $L^2$-index theorem. Ma and Marinescu \cite[3.6]{MM07} reduced the covering case to the compact case by comparing heat kernels, however our following approach is direct from Theorem \ref{heatkernelthm}. 

We firstly recall the notations for coverings.
Let $(\widetilde{M}, J)$ be a complex manifold of dimension $n$ with a compatible Riemannian metric $g^{T\widetilde{M}}$. Let $\widetilde \omega$ be the associated real $(1,1)$-form defined by $\widetilde \omega(u,v)=g^{T\widetilde{M}}(J u,v)$ on ${T\widetilde{M}}$.
A group $\Gamma$ is called a discrete group acting holomorphically, freely and properly on $\widetilde{M}$,
if $\Gamma$ is equipped with the discrete topology such that (1) the map  $\Gamma \times \widetilde{M} \rightarrow \widetilde{M}, (r, \widetilde{x})\mapsto r.\widetilde{x}$ is holomorphic; (2) $r.\widetilde{x}=\widetilde{x}$ for some $\widetilde{x} \in \widetilde{M}$ implies that $r=e$ the unit element of $\Gamma$; and (3) the map  $\Gamma \times \widetilde{M} \rightarrow \widetilde{M}\times \widetilde{M}, (r,\widetilde{x})\rightarrow (r.\widetilde{x},\widetilde{x})$ is proper.

A Riemannian metric $g^{T\widetilde{M}}$ (or $\widetilde{\omega}$) is called $\Gamma$-equivariant, if the map $r:\widetilde{M}\rightarrow \widetilde{M}$ is isometric with respect to $g^{T\widetilde{M}}$
for every $r\in \Gamma$. We say a Hermitian manifold $(\widetilde{M},\widetilde{\omega})$ is a covering manifold, if there exists a discrete group $\Gamma$ acting holomorpically, freely and properly on $\widetilde{M}$ such that $\widetilde{\omega}$ is $\Gamma$-equivariant and the quotient $M:=\widetilde{M}/\Gamma$ is compact. 

A holomorphic Hermitian vector bundle $(\widetilde{F},h^{\widetilde{F}})$ over $\widetilde{M}$ is called $\Gamma$-equivariant, if there is a map $r_{\widetilde{F}}:{\widetilde{F}}\rightarrow {\widetilde{F}}$ associated to $r\in \Gamma$, which commutes with the fibre projection $\pi:{\widetilde{F}}\rightarrow \widetilde{M}$ (i.e., $r\circ\pi =\pi\circ r_{\widetilde{F}}$), such that $h^{\widetilde{F}}(v,w)=h^{\widetilde{F}}(r_{{\widetilde{F}}}v,r_{{\widetilde{F}}}w)$ for $v,w \in \widetilde{F}$. 

A relatively compact open subset $U \Subset \widetilde{M} $ is called a fundamental domain of the action $\Gamma$ on a covering manifold $ \widetilde{M}$, if the following conditions are satisfied: (a) $ \widetilde{M}=\cup_{r\in\Gamma}r(\overline{U}) $; (b) $r_1(U)\cap r_2(U)$ is empty for $ r_1,r_2 \in \Gamma$ with $r_1\neq r_2$; and (c) $\overline{U}\setminus U$ has zero measure. The fundamental domain always exists and there is a biholomorphic map between $U$ and $M\setminus Z$ with a zero measure subset $Z\subset M$.

\begin{thm}[\cite{TCM:01,MTC:02,MM07}]\label{thm_cplxcover}
Let $(\wt M,\widetilde \omega)$ be a Hermitian manifold of dimension $n$. Let $\Gamma$ be a discrete group acting  holomorphically, freely, and properly on $\widetilde M$, such that $\widetilde{\omega}$ is a $\Gamma$-equivariant Hermitian
metric and the quotient $M=\widetilde M/\Gamma$ is compact. 
Let $(\widetilde L,h^{\widetilde L})$ be a $\Gamma$-equivariant
	holomorphic Hermitian line bundle on $\widetilde M$. Let $0\leq q\leq n$. Then, as $k\rightarrow \infty$, the strong holomorphic Morse inequalities hold for the von Neumann dimension of the reduced $L^2$-Dolbeault cohomology,
	\begin{equation}\nonumber%\label{eq_hmi_cover}
	\sum_{j=0}^q(-1)^{q-j}\dim_{\Gamma}\overline{H}^{j}_{(2)}(\widetilde{M},\widetilde{L}^k)\leq \frac{k^n}{n!}\int_{M(\leq q)}(-1)^q c_1(L,h^L)^n+o(k^n),
	\end{equation}
	with equality for $q=n$, and the weak Morse inequalities hold
	\begin{equation}\nonumber
	\begin{split}
	\dim_\Gamma \ov H^{q}_{(2)}(\til M,\til L^k)\leq \frac{k^n}{n!}\int_{M(q)}(-1)^qc_1(L,h^L)^n+o(k^n),
	\end{split}
	\end{equation}  
	where $M(q)$ is the subset of $M$ consisting of points on which $R^L$ is non-degenerate and has exactly $q$ negative eigenvalues. 
\end{thm} 

\begin{proof}
	We apply Theorem \ref{heatkernelthm} on the fundamental domain $U\Subset \til M$ as \eqref{eq-trlocal}, then we use the general fact \cite[Lemma 3.6.6,(3.6.23)]{MM07}: for any $t>0$, $k\in\N^*$, $0\leq q\leq n$,  
	\begin{equation} \nonumber
	\sum_{j=0}^q(-1)^{q-j}\dim_{\Gamma}\overline{H}_{(2)}^j(\til M,\til L^k)\leq \sum_{j=0}^q(-1)^{q-j}\Tr_{\Gamma}\Big(e^{-\frac{t}{k}\Box_{k}^j}\Big)=\sum_{j=0}^q(-1)^{q-j}\int_U \Tr\Big( e^{-\frac{t}{k}\Box_{k}^j}(z,z)\Big)dv_{\til M}(z),
	\end{equation}  
	with equality for $q=n$, where $\overline{H}_{(2)}^j(\til M,\til L^k)$ is the reduced $L^2$ Dolbeault cohomology.
\end{proof}

 The goal of this section is to establish holomorphic Morse inequalities for symmetric tensor power of vector bundle on covering manifolds. For this purpose, we next establish a $L^2$ version of Le Potier isomorphism for holomorphic sections of symmetric tensor power of holomorphic vector bundles.
  Let $M$ be a complex manifold of dimension $n$.  Let $E$ be a holomorphic vector bundle of rank $r$ on $M$. Let $E^*$ be the dual bundle of $E$. We denote by $P(E^*)$ the projective bundle associated to $E^*$, which is a complex manifold of dimension $n+r-1$. Let $\cO_{P(E^*)}(-1)$ be the tautological line bundle on $P(E^*)$ and $\cO_{P(E^*)}(1)$  its dual bundle. Let $p\in \N
  \setminus\{0\}$ and let $S^p(E)$ be the $p$-th symmetric tensor power of $E$ and $\cO_{P(E^*)}(p)$ the $p$-th tensor power of $\cO_{P(E^*)}(1)$. The theorem of Le Potier (\cite[Chap.III \S 5 (5.7)]{Kob:87}) implies that there exists an isomorphism as follows
  \begin{equation}\label{aeq1}
  	\rho: H^{0} (M, S^p(E))\simeq H^{0} (P(E^*), \cO_{P(E^*)}(p)), \quad S\mapsto \rho(S)=\bs.
  \end{equation}

  For $x\in M$, we denote by $E_x$ the fibre of $E$ at $x$. Let $h^E$ be Hermitian metric on $E$. For $v\in E_x$, its metric dual $v^*$ is given by $v^*=h^{E}(\cdot, v)$. And we denote the space of such dual vectors by $E_x^*$. Let $S^p(E^*_x)$ be the $p$-th symmetric tensor power of $E^*_x$ and let $P(E^*_x)$ be the projectlization. For $v\in E_x\setminus \{0\}$, we denote the equivalent class of  $v^*$ by $[v^*]\in P(E^*_x)\cong P(\C^r)$ under an orthonormal basis of $(E_x,h^{E_x})$, and it is clear that $v^{*\otimes p}\in S^p(E^*_x)$. From (\ref{aeq1}), we have
  $$\bs([v^*])(v^{*\otimes p})=S(x)(v^{*\otimes p})\in \C.$$
  	Let $h^{S^p(E)}$ and $h^{\cO_{P(E^*)}(-p)}$ be Hermitian metrics on the induced bundles by $(E,h^E)$ respectively.
  It follows that, for $v\in E_x\setminus \{0\}$ with $|v|_{h^E}=1$ and $e_v\in \cO_{E^*}(1)|_{[v^*]}$ such that $e_v(v^*)=1$,
  $$\bs([v^*])=\langle S(x),v^{\otimes p}\rangle_{h^{S^p(E)}} e_v^{\otimes p}.$$
  Moreover, from $\langle e_v,e_v \rangle_{h^{\cO_{P(E^*)}(1)}}=1$, we have for any $S, T \in H^0(M, S^p(E))$,
  	\begin{equation}\label{aeq2}
  		\langle \bs([v^*]),\bt([v^*]) \rangle_{h^{\cO_{P(E^*)}(p)}}=\langle S(x), v^{\otimes p} \rangle_{h^{S^p(E)}} \langle v^{\otimes p}, T(x) \rangle_{h^{S^p(E)}}.
  	\end{equation}

  Let $\pi:P(E^*)\rightarrow M$ be the natural projection with fibers $\pi^{-1}(x)=P(E_x^*)$. Let $\omega$ be Hermitian metric on $M$. The prescribed Hermitian metric on $P(E^*)$ is defined by
  \begin{equation}\label{eq_prevol}
  	\omega_{P(E^*)}:=\omega(x)+\omega_{P(E^*_x)}([v^*]).
  \end{equation} 
  The induced volume form is given by
  $dV_{P(E^*)}(x,[v^*])=dV_M(x)\wedge dV_{P(E^*_x)}([v^*])$, where $dV_M=\omega^n/n!$ and $dV_{P(E^*_x)}=\omega_{P(E^*_x)}^{r-1}/(r-1)!$. 
  
  \begin{prop}\label{eq_l2lep} With respect to the metric  \eqref{eq_prevol}, we have the $L^2$ version of Le Potier isomorphism 
      \begin{eqnarray}\nonumber
  	H_{(2)}^{0} (M, S^p(E))&\simeq& H_{(2)}^{0} (P(E^*), \cO_{P(E^*)}(p)),\quad p\in\N\setminus\{0\}. 
  \end{eqnarray}
  \end{prop}
In fact, Proposition \ref{eq_l2lep} reduces to the following observation.
  \begin{thm}\label{aprop1}
  	Let $(M,\omega)$ be a Hermitian manifold of dimension $n$.  Let $(E,h^E)$ be a holomorphic Hermitian vector bundle of rank $r$ on $M$.
  	Then, for any $S,T\in H^{0} (M, S^p(E))$ for $p\in \N\setminus\{0\}$, we have 
  	\begin{equation}\label{aeq3}
  		 \left(\bs,\bt\right)_{L^2 (M, S^p(E))}=\frac{p!}{(p+r-1)!}\left(S, T\right)_{L^2 (P(E^*), \cO_{P(E^*)}(p))},
  	\end{equation}
  	where the $L^2$ inner products are with respect to the above volume forms $dV_M$ and $dV_{P(E^*)}$.  
    %induced by the metric $\omega$ and \eqref{eq_prevol} respectively.
  \end{thm}
  \begin{proof} 
    For simplifying notations, we write the Hermitian metrics by $\langle\,\cdot\, , \,\cdot\,\rangle_h$ and the $L^2$ inner products by $(\,\cdot\, , \,\cdot\,)_{L^2}$.
  	Given $x\in M$, we assume $S(x)\neq 0$. Set 
  	$\lambda(x):=|S(x)|_h>0.$ 
   Without loss of generality (\cite[Lemma 4.2]{CGLM:08}) we can assume 
   $S(x)=\lambda (x)e_1^{\otimes p},$
  	where $e_1\in E_x$ with unit norm.     
  We write
  	\begin{equation}\label{aeq4}
  		(\bs, \bt)_{L^2}
  		=\int_M \fc(x) dV_M,
  	\end{equation}
  	where $\fc(x):=\int_{P(E_x^*)}\langle \bs([v^*]), \bt([v^*]) \rangle_h dV_{P(E_x^*)}$ is the integration along the fiber at $x\in M$. By (\ref{aeq2}), 
  	\begin{equation}\label{aeq5}
  		\fc(x)=\int_{|v|_{h=1}, [v^*]\in P(E_x^*)} \langle S(x), v^{\otimes p} \rangle_h \langle v^{\otimes p}, T(x) \rangle_h dV_{P(E_x^*)}.
  	\end{equation}
  	
  	We extend $e_1$ to an orthonormal basis $\{ e_i \}_{i=1}^r$ of $E_x$. For any unit norm vector $v=\sum_1^r v^ie_i\in E_x$,
  	\begin{equation}\label{aeq6}
  		\langle S(x), v^{\otimes p} \rangle_h=\lambda(x)\overline{v^1}^{p}.
  	\end{equation}
  	Moreover, we write
  	$T(x)=\sum_k b_k^1\otimes ...\otimes b_k^p \in S^p(E_x),$ where $b_k^l=\sum_{i=1}^r b_k^{l,i}e_i\in E_x$.
  	Then 
  	\begin{equation}\label{aeq7}
  		\langle v^{\otimes p}, T(x) \rangle_h=\sum_k\prod_{l=1}^p (\sum_{i=1}^r v^i \overline{b_k^{l,i}}).
  	\end{equation}
  	 
  	Next we consider $[v^*]\in P(E^*_x)$ with 
  	$v^*=\sum_{i=1}^r\overline{v^i}e_i^*$ such that $\sum_{i=1}^r |v^i|^2=1$. By (\ref{aeq6}) we assume $\overline{v^1}\neq 0$, the integral area of  $P(E_x^*)$ in (\ref{aeq5}) consists of $(v^1,u^2,...,u^r)\in \C^r$ as follows: 
  	$$\overline{v^1}\neq 0,\quad u^j:=\frac{\overline{v^j}}{\overline{v^1}},~j=2,...,r$$ 
  	satisfying
  	$1+\sum_{j=2}^{r}|u^j|^2=|v^1|^{-2}.$
  	Set $x_j+\sqrt{-1}y_j=u^{j+1}$ for $j=1,...,r-1$, it follows that 
  	\begin{equation}\label{aeq8}
  		dV_{P(E_x^*)}=\frac{\omega_{FS}^{r-1}}{(r-1)!}
  		=\frac{dx_1\wedge dy_1\wedge...\wedge dx_{r-1}\wedge dy_{r-1}}{\pi^{r-1}(1+\sum_1^{r-1} (|x_j|^2+|y_j|^2) )^r}. 
  	\end{equation}
  	
  	By combining (\ref{aeq6}) and (\ref{aeq7}), we have 
  	\begin{equation}\label{aeq9}
  		\langle S(x), v^{\otimes p} \rangle_h\langle v^{\otimes p}, T(x) \rangle_h
  		=\lambda(x)\sum_k \frac{\prod_{l=1}^p (\overline{b_k^{l,1}}+ \sum_{t=1}^{r-1} ( x_t-\sqrt{-1}y_t ) \overline{ b_k^{l,t+1} })}
  		{ ( 1+\sum_{t=1}^{r-1}( |x_t|^2+|y_t|^2 )   )^p   } . 
  	\end{equation}
  	We substitute (\ref{aeq8}) and (\ref{aeq9}) into (\ref{aeq5}).
  	Moreover, we change the coordinates of $\R^{2r-2}$ by
  	$\eta_te^{i\theta_t}=x_t+\sqrt{-1}y_t$ for $1\leq t\leq r-1$ and obtain
  	\begin{equation}\nonumber%\label{aeq11}
  		\fc(x)=\frac{\lambda(x)}{\pi^{r-1}} \sum_k \int_{\R^{2r-2}} 
  		\frac{\prod_{l=1}^p (\overline{b_k^{l,1}}+ \sum_{t=1}^{r-1} ( \eta_te^{-i\theta_t}) \overline{ b_k^{l,t+1} })}
  		{ ( 1+\sum_{t=1}^{r-1}\eta_t^2   )^{p+r} } \eta_1d\eta_1\wedge d\theta_1\wedge...\wedge\eta_{r-1}d\eta_{r-1}\wedge d\theta_{r-1}.
  	\end{equation}
  	By the fact that 
  	$\int_{0\leq \theta\leq 2\pi} e^{-i\theta}d\theta=0$,
  	it follows that 
  	\begin{equation}\nonumber%\label{aeq12}
   \begin{split}
  		\fc(x)
  		&= \frac{\lambda(x)}{\pi^{r-1}} (\sum_k\prod_{l=1}^p \overline{b_k^{l,1}}) (2\pi)^{r-1} \int_{(\R^+)^{r-1}}
  		\frac{\eta_1...\eta_{r-1}}
  		{ ( 1+\sum_{t=1}^{r-1}\eta_t^2   )^{p+r} } d\eta_1\wedge...\wedge d\eta_{r-1} \\
  		&=\lambda(x)\langle e_1^{\otimes p}, T(x) \rangle_h 2^{r-1}  \int_{(\R^+)^{r-1}} 
  		\frac{\eta_1...\eta_{r-1}}
  		{ ( 1+\sum_{t=1}^{r-1}\eta_t^2   )^{p+r} } d\eta_1\wedge...\wedge d\eta_{r-1} \\
  		&= \frac{\langle S(x), T(x) \rangle_h }{(p+r-1)(p+r-2)...(p+1)}.   
   \end{split}
  \end{equation}
  Finally, we substitute it into (\ref{aeq4}), and then (\ref{aeq3}) follows.
  \end{proof}
\noindent\textbf{Proof of Theorem \ref{cor_sym}}: It follows from Proposition \ref{eq_l2lep} and Theorem \ref{thm_cplxcover} for ${\cO_{\til{P(E^*)}}(k)}$ over $\til{P(E^*)}$.\qed
  
%\begin{cor}
%	If $\int_{P(E^*)(\leq 1)} c_1(\cO_{P(E^*)}(1),h^{\cO_{P(E^*)}(1)})^{n+r-1}>0$, then the $L^2$-volume (See Def. \ref{def-gamvol})
%	$$\vol_{\Gamma,(2)}(\til E,h^{\til E})>0.$$  
%\end{cor}     
 
 % \cite[Theorem 3.2.16,Theorem 3.5.10(3.5.25)(3.5.26)]{MM07}. 

	Based on \cite[Theorem 3.3.5(ii), Theorem 3.5.12, Corollary 3.5.13]{MM07}, the Le Potier isomorphism \eqref{aeq1} and its $L^2$ version (Proposition \ref{eq_l2lep}) also lead to holomorphic Morse inequalities for symmetric tensor power of vector bundle on non-compact complex manifolds as follows.
  
  \begin{thm}
  Let $(M,\omega)$ be a Hermitian manifold of dimension $n$ and $(E,h^E)$ be a holomorphic Hermitian vector bundle on $M$ of rank $r$. Suppose $\omega_{P(E^*)}=\sqrt{-1}R^{\cO_{P(E^*)}(1)}>0$ defines a complete Hermitian metric on $P(E^*)$ on which $\mbox{Ric}_{\omega_{P(E^*)}}\geq -C\omega_{P(E^*)}$. Then, as $k\rightarrow\infty$, we have
  \begin{equation}\nonumber
  	\dim H_{(2)}^{0} (M, S^k(E))\geq  \frac{k^{(n+r-1)}}{(n+r-1)!}\int_{P(E^*)(\leq 1)} c_1(\cO_{P(E^*)}(1),h^{\cO_{P(E^*)}(1)})^{n+r-1}+o(k^{n+r-1}).   
  \end{equation}
\end{thm}

\begin{thm}
  Let $M$ be a complex manifold of dimension $n$ and $(E,h^E)$ be a holomorphic Hermitian vector bundle on $M$ of rank $r$. Suppose $P(E^*)$ be a weakly $1$-complete manifold.
  Assume  $(\cO_{P(E^*)}(1),h^{\cO_{P(E^*)}(1)})$ is positive outside a compact set $K\subset P(E^*)$. Then, as $k\rightarrow\infty$, we have
  
  \begin{equation}\nonumber
  	\begin{split}
  &	\sum_{j=r}^{n+r-1}(-1)^{j-r}\dim H^j(M,S^k(E))\\
   &\leq\frac{k^{(n+r-1)}}{(n+r-1)!}\int_{P(E^*)(\geq r)}(-1)^r c_1(\cO_{P(E^*)}(1),h^{\cO_{P(E^*)}(1)})^{n+r-1}+o(k^{n+r-1}),\\
  	& \dim H^0(M,S^k(E))\geq\frac{k^{(n+r-1)}}{(n+r-1)!}\int_{P(E^*)(\leq 1)}c_1(\cO_{P(E^*)}(1),h^{\cO_{P(E^*)}(1)})^{n+r-1}+o(k^{n+r-1}).
  	\end{split}
  \end{equation}
In particular, if $(\cO_{P(E^*)}(1),h^{\cO_{P(E^*)}(1)})$ is positive, then
  %\begin{equation}\nonumber
  	$\lim_{k\rightarrow\infty}k^{-(n+r-1)}\dim H^0(M,S^k(E))=+\infty.$
  %\end{equation} 
\end{thm}

 \subsection{Algebraic Morse inequalities on K\"{a}hler coverings}\label{sec_algmi}
 		\
	
\vspace{0.18cm}
% \subsection{Volume of nef line bundles with Lie group action}
 The purpose of this section is to generalize the formalism on the volume of nef line bundles (see Corollary \ref{thm_cover}) and algebraic Morse inequalities involved only algebraic invariants to K\"{a}hler coverings (see Theorem \ref{thm_amikc}) via Theorem \ref{thm_cplxcover}. The projective and compact K\"{a}hler manifold cases were due to Trapani, Siu, Angelini and Demailly. 
   
  Let $(M,\omega)$ be a Hermitian manifold of dimension $n$. Let $1_A$ be the characteristic function for a subset $A\subset M$, i.e., $1_A(x)=1$ for $x\in A$; $1_A(x)=0$ for $x\in M\setminus A$. We next collect two inequalities essentially due to Demailly and Siu. The proof follows from the linear algebraic computation in an orthonormal frame of $T^{1,0}M$ and the induction process on the number $n$.   
\begin{itemize}
     \item[(i)]  Let $\epsilon>0$ and $\Theta_\epsilon$ be a real $(1,1)$-from on $M$. Suppose $\Theta_\epsilon\geq -\epsilon \omega$ for some Hermitian form $\omega$ on $M$. Then,
 	at every point of $M$ we have
 	\begin{equation}\label{eq_1}
 		0\leq (-1)^q 1_{M(q,\Theta_\epsilon)}\frac{\Theta_\epsilon^n}{n!}  \leq \frac{(\epsilon\omega)^q}{q!}\wedge \frac{(\Theta_\epsilon+\epsilon\omega)^{n-q}}{(n-q)!},
 	\end{equation}
 where $M(q,\Theta_\epsilon)$ is the set of point $x$ of $M$ such that $\Theta_\epsilon$ is non-degenerated and has exactly $q$ negative eigenvalues at $x$.
 \item[(ii)] Let $\lambda_j$ be positive smooth functions on $M$ for $j=1,\cdots,n$. Let $\lambda_1\geq \lambda_2\geq \cdots\geq \lambda_n>0$ on $M$. Let
 	 $f_n^j(\lambda)$ be the $j$-th elementary symmetric polynomial in $\lambda_1,\cdots,\lambda_n$ given by
 	%\begin{equation}\nonumber
 		$f_n^j(\lambda):=\sum_{1\leq l_1<\cdots<l_j\leq n}\lambda_{l_1}\cdots\lambda_{l_j}$,
 	%\end{equation}
 	where $f^0_n(\lambda)=1$.
 	If $0\leq  q\leq n$, then at each point of $M$ we have  
 	\begin{equation}\label{eq_2}
 		1_{\{x\in M: \lambda_{q+1}(x)<1 \}}(-1)^q\prod_{j=1}^n(1-\lambda_j)	\leq \sum_{j=0}^q(-1)^{q-j}f^j_n(\lambda),
 	\end{equation}
  with equality for $q=n$.  
 \end{itemize}

 	A $\Gamma$-equivariant holomorphic line bundle $\til L$ is called $\Gamma$-nef, if for a $\Gamma$-equivariant Hermitian metric $\til \omega$ on $\til M$, for every $\epsilon>0$ there exists a smooth $\Gamma$-equivariant Hermitian metric $h^{\til L}_\epsilon$ on $\til L$ such that  $c_1(L,h^{L}_\epsilon)\geq -\epsilon\omega$ on the quotient $M=\til M/\Gamma$. 	
 When $\Gamma$ is trivial, it reduces to the usual definition of nef line bundle. For simplifying, we denote $\limsup_{\epsilon\rightarrow 0,p\rightarrow \infty}:=\limsup_{\epsilon\rightarrow 0}\limsup_{p\rightarrow \infty}$ in the following.
 \begin{lem}\label{lem_gam_nef}
 	If $(\til M,\til \omega)$ is a K\"{a}hler covering manifold of $\dim \til M=n$ and $\til L$ is a $\Gamma$-nef line bundle, then
 	\begin{equation}\nonumber
 	\limsup_{\epsilon\rightarrow 0,p\rightarrow \infty} p^{-n}\dim_{\Gamma} \ov H_{(2)}^{q}(\til M,\til L^p, h^{\til L^p}_\epsilon)=0, \quad 1\leq q\leq n.
 	\end{equation}	
 \end{lem} 
 \begin{proof}
 	Let $\pi_\Gamma: \til M\rightarrow M$ be the natural projection.	For each $\epsilon>0$, there exits a metric $h^{\til L}_\epsilon$ such that $c_1(L,h^{L}_\epsilon)\geq -\epsilon \omega$ on $M$, where $\pi_\Gamma^* (\omega)=\til \omega$ and $\pi^*_\Gamma(c_1( L,h^{L}_\epsilon))=c_1(\til L,h^{\til L}_\epsilon)$.
 	We set $\Theta_\epsilon:=c_1( L,h^{ L}_\epsilon)$ and
 	apply Theorem \ref{thm_cplxcover}, it follows that
 	\begin{equation}\nonumber
 	\limsup_{p\rightarrow\infty} p^{-n}\dim_{\Gamma} \ov H_{(2)}^{q}(\til M,\til L^p,h^{\til L^p}_\epsilon)\leq \frac{1}{n!}\int_{M(q,\Theta_\epsilon))}(-1)^q \Theta_\epsilon^n.
 	\end{equation}
 	%\end{equation}
 	For a $d$-closed $2$-form $\alpha$ on $M$ we denote by $[\alpha]$ the cohomology class. Note the class $[\Theta_\epsilon]=c_1(L)$ is independent of $\epsilon$. Since $d\omega=0$, the compactness of $M$ and \eqref{eq_1}, we obtain
 	\begin{equation}\nonumber
 	\limsup_{p\rightarrow\infty} p^{-n}\dim_{\Gamma} \ov H_{(2)}^{q}(\til M,\til L^p,h^{\til L^p}_\epsilon)\leq
 	\frac{\epsilon^q}{q!(n-q)!}\int_{M}[\omega]^q\wedge ([c_1(L)]+\epsilon[\omega])^{n-q}.
 	\end{equation}
 	 for $1\leq  q\leq n$. The proof is complete by letting $\epsilon\rightarrow 0$.  
 \end{proof}  

 The volume of a $\Gamma$-equivariant holomorphic Hermitian line bundle $(\til L,h^{\til L})$ is defined by
 	$$\vol_{\Gamma, (2)}(\til L,h^{\til L}):=\lim\sup_{p\rightarrow \infty}\frac{n!}{p^n}\dim_\Gamma H^0_{(2)}(\til M,\til L^p,h^{\til L^p}).$$
Since the von Neumann dimension can be viewed as a normalized dimension, $\vol_{\Gamma,(2)}(\til L,h^{\til L})$ can be viewed as a normalized volume. In particular, when $\Gamma$ is finite,
 $\vol_{\Gamma,(2)}(\til L,h^{\til L})=|\Gamma|^{-1}\vol(\til L)$; when $\Gamma$ is trivial, it reduces to the usual volume of line bundle. Furthermore, the volume of a $\Gamma$-nef line bundle $\til L$ on K\"{a}hler coverings can be defined by
\begin{equation}\label{def-gamnef} 
	\vol_{\Gamma,(2)}(\til L):=\lim_{\epsilon\rightarrow 0}\vol_{\Gamma, (2)}(\til L,h^{\til L}_\epsilon).
\end{equation}
In fact, from Lemma \ref{lem_gam_nef} and Atiyah's $L^2$-index theorem, it follows that as $\epsilon\rightarrow 0$,
 \begin{equation}\label{eq-weldefvol}
 \vol_{\Gamma,(2)}(\til L,h^{\til   L}_\epsilon)=\limsup_{p\rightarrow\infty}\frac{n!}{p^n} \dim_{\Gamma}H^0_{(2)}(\til M,\til L^p,h^{\til L^p}_\epsilon)=\int_{M}c_1(L)^n+o(\epsilon).
 \end{equation}

\noindent\textbf{Proof of Theorem \ref{thm_cover}}: 
The first assertion follows from \eqref{def-gamnef} and \eqref{eq-weldefvol}. The second assertion holds for K\"{a}hler forms $c_1(L_j,h^{L_j}_\epsilon)+\epsilon \omega>0$ with $j=1,2$ and let $\epsilon\rightarrow 0$ as in the proof of Lemma \ref{lem_gam_nef}.\qed

Analogous to the compact case \cite[(2.3) Theorem]{Dem10}, we generalize the algebraic Morse inequalities to K\"{a}hler coverings via \eqref{eq_2} and Theorem \ref{thm_cplxcover}.  
  
 \begin{thm}[Algebraic Morse inequalities on K\"{a}hler coverings]\label{thm_amikc}
 	Let $\til M$ be a K\"{a}hler covering manifold of $\dim \til M=n$. Let $\til F$ and $\til E$ be $\Gamma$-nef line bundles on $\til M$.  	
 	Then we have 
 	\begin{itemize}
 		\item [(a)] Strong algebraic Morse inequalities hold: For any $q=0,1,\cdots,n$, 
 		\begin{equation}\nonumber
 		\begin{split}
 		\limsup_{\epsilon\rightarrow 0,p\rightarrow \infty}n!p^{-n}\chi^q_{\Gamma} (\til M, (\til F\otimes \til E^*)^p, h^{(\til F\otimes \til E^*)^p}_\epsilon )
 		\leq \sum_{j=0}^{q}(-1)^{q-j}\binom{n}{j}\int_M c_1(F)^{n-j}\wedge c_1(E)^j
 		\end{split}
 		\end{equation}
 		with equality for $q=n$, where $\chi^q_{\Gamma}(\til M, \til E,h^{\til E}):=\sum_{j=0}^{q}(-1)^{q-j}\dim_{\Gamma} H^{j}_{(2)}(\til M, \til E,h^{\til E})$.
 		
 		\item [(b)] Weak algebraic Morse inequalities hold: For any $q=0,1,\cdots,n$, 
 		\begin{equation}\nonumber
 		\begin{split}
 		\limsup_{\epsilon\rightarrow 0,p\rightarrow \infty}n!p^{-n}\dim_{\Gamma} H^{q}_{(2)}(\til M, (\til F\otimes \til E^*)^p,h^{(\til F\otimes \til E^*)^p}_\epsilon )
 		\leq \binom{n}{q}\int_M c_1(F)^{n-q}\wedge c_1(E)^q.
 		\end{split}
 		\end{equation}
 		
 		\item [(c)] The volume of the $\Gamma$-nef line bundle $\til F\otimes\til E^*$:
 		\begin{equation}\nonumber
 		\begin{split}
 		\int_M c_1(F)^n\geq \vol_{\Gamma,(2)}(\til F\otimes\til E^*) 
 		&\geq \int_M c_1(F)^n-nc_1(F)^{n-1}\wedge c_1(E).
 		\end{split}
 		\end{equation}
 	\end{itemize} 
 \end{thm}

 \begin{cor}   
 	Let $\til M$ be K\"{a}hler covering manifold of $\dim \til M=n$. Let $\til F$ and $\til E$ be $\Gamma$-nef line bundles on $\til M$.  Suppose $F$ is big, i.e., $\int_M c_1(F)^n>0$.
 	Then some tensor power of $\til F^k\otimes \til E^*$ has a non-trivial $L^2$-holomorphic section on $\til M$, when
 	\begin{equation}\nonumber
 	k>n\Big(\int_M c_1(F)^{n-1}\wedge c_1(E) \Big)\Big(\int_M c_1(F)^n \Big)^{-1}.
 	\end{equation} 
 \end{cor}   
 
  \subsection{CR covering manifolds}
	In the same vein as our proof of Theorem \ref{thm_cplxcover}, by applying heat kernel asymptotics for Kohn Laplacian \cite[Theorem 1.1]{HZ23} to CR covering manifolds, we obtain new Morse inequalities as a generalization of Morse inequalities for compact CR manifolds \cite{HM12} (See also \cite{HS20} for the $S^1$-action case).
	\begin{thm}\label{thm_crcover}
		Let $\til X$ be an orientable CR manifold of dimension $2n+1$, $n\geq 1$, on which a discrete group $\Gamma$ acts smoothly, freely, and properly such that the quotient $X=\til X/\Gamma$ is a compact CR manifold. Let $\til \omega$ be a $\Gamma$-equivariant Hermitian metric on $\til X$. Let $(\widetilde L,h^{\widetilde L})$ be a $\Gamma$-equivariant
		CR Hermitian line bundle on $\widetilde X$.
		Let $0\leq q\leq n$ and suppose that condition $Y(j)$ holds at each point of $\til X$ for all $0\leq j\leq q$. Then, as $k\rightarrow \infty$, the strong CR Morse inequalities hold for the von Neumann dimension of the reduced $L^2$-Kohn-Rossi cohomology,
		\begin{equation}\nonumber
			\sum_{j=0}^q(-1)^{q-j}\dim_{\Gamma}\overline{H}_{b,(2)}^j(\til X,\til L^k)\leq \frac{k^{n+1}}{(2\pi)^{n+1}}\int_X \int_{\R_x(\leq q)}(-1)^{q}\det(\dot{\mathcal R^\phi_x}-2t\dot{\mathcal L_x})dtdv_X(x)+o(k^{n+1})
		\end{equation}  
	with equality for $q=n$. In particular, for $q=1$, we have
	\begin{equation}\nonumber 
	\lim_{k\to+\infty}k^{-(n+1)}\dim_{\Gamma}\overline{H}_{b,(2)}^0(\til X,\til L^k)\geq \frac{1}{(2\pi)^{n+1}}\int_X \int_{\R_x(\leq 1)}\det(\dot{\mathcal R^\phi_x}-2t\dot{\mathcal L_x})dtdv_X(x),
	\end{equation}
		where $\R_x(j)=\{t\in\R:\dot{\mathcal R^\phi_x}-2t\dot{\mathcal L_x}~\mbox{is~non-degenerated and has}~j~\mbox{negative eigenvalues}\}$ and $\R_x(\leq q)=\cup_{0\leq j\leq q}\R_x(j)$.
	\end{thm}   

	%The Morse inequalities on CR coverings with the $S^1$-action were studied in \cite{HS20}, and Shuang Su also considered heat kernel asymptotics of Kodaira Laplacian in a different way. It is interesting to extend above result to CR vector bundle case.

\end{document}